\documentclass{article}

\usepackage{hyperref}

\usepackage{amsfonts}
\usepackage{amssymb}
\usepackage{textcomp}
\usepackage{latexsym,amsmath}
\usepackage{graphics}

\renewcommand{\phi}{\varphi}
\newcommand{\be}{\begin{equation}}
\newcommand{\ee}{\end{equation}}
\newcommand{\ba}{\begin{eqnarray}}
\newcommand{\ea}{\end{eqnarray}}
\newcommand{\ban}{\begin{eqnarray*}}
\newcommand{\ean}{\end{eqnarray*}}

\newcommand{\nul}{{\bf0}}
\newcommand{\rd}{{\mathbb R}^d}
\newcommand{\zd}{{\mathbb Z}^d}
\newcommand{\td}{{\mathbb T}^d}

\newcommand{\z} {{\mathbb Z}}

\newcommand{\n} {{\mathbb N}}

\newcommand{\h}{\widehat}
\newcommand{\w}{\widetilde}

\def\spec{\operatorname{spec}}

\def\N{{{\mathbb N}}}
\def\Z{{{\mathbb Z}}}
\def\T{{{\mathbb T}}}
\def\R{{\mathbb R}}

\def\vp{{\varphi}}

\def\({\left(}
\def\){\right)}

\newtheorem{theo}{Theorem}
\newtheorem{lem}[theo]{Lemma}
\newtheorem {prop} [theo] {Proposition}
\newtheorem {coro} [theo] {Corollary}
\newtheorem {defi} [theo] {Definition}
\newtheorem {rem} [theo] {Remark}

\title{Approximation by periodic multivariate \\ quasi-projection operators
\thanks{The first author was partially supported by DFG project KO 5804/1-1 (Section 4.2 belongs to this author); the second and the third  authors are supported by the Russian Science Foundation under grant No. 18-11-00055 (Sections~4.1 and 5 belong to these authors) }}
\author{
Yu. Kolomoitsev$^{1, 2}$, A. Krivoshein$^{3}$ and M. Skopina$^{3}$
}
\date{\small $^{1}$Universit\"at zu L\"ubeck, Institut f\"ur Mathematik, L\"ubeck, Germany \\
\small $^{2}$Institute of Applied Mathematics and Mechanics of NAS of Ukraine, Slov'yans'k, Ukraine\\
 $^{3}$St. Petersburg State University, Russia \\
kolomoitsev@math.uni-luebeck.de, a.krivoshein@spbu.ru, skopina@ms1167.spb.edu}

\begin{document}

\maketitle

\begin{abstract}
Approximation properties of periodic quasi-projection operators with matrix dilations are studied. Such operators are generated by a sequence of functions  $\phi_j$ and a sequence of distributions/functions $\w\phi_j$.  Error estimates for sampling-type quasi-projection operators are obtained  under the periodic  Strang-Fix conditions for $\phi_j$ and the compatibility conditions  for  $\phi_j$ and  $\w\phi_j$.
These estimates are given in terms of the Fourier coefficients of approximated functions and provide  analogs of some known non-periodic  results. Under some additional assumptions error estimates are given in other terms in particular using the best approximation.
A number of examples  are  provided.
\end{abstract}

\bigskip


\medskip


\section{Introduction}

Approximation properties of non-periodic quasi-projection operators are actively studied
 by many authors. The class of such operators is very
large, it includes classical sampling and sampling-type expansions  (see, e.g.,~\cite{Unser, Zayed, Brown, Si2,  Butz4, Butz5, ButzB, JZ,  KS, KKS, Lif} and the references therein), Kantorovich-Kotelnikov operators and their generalizations (see, e.g.,~\cite{ Butz6,  OT15,CV0, CV1, CV2, VZ2, KS3}), scaling expansions associated with wavelet constructions  (see, e.g.,~\cite{Jia1, Jia2,  BDR, DB-DV1, v58, KS0, Sk1}) and others.
The most general form of the multivariate quasi-projection operator with a dilation matrix  $M$ is given by
\be
 \sum_{k\in\zd} m^{j} \langle f,\w\vp(M^j\cdot-k)\rangle \vp(M^j\cdot-k),
\label{my00}
\ee
where $\phi$ is a function and $\w\phi$ is a distribution or function,  $m=|\det M|$, and the inner product $\langle f,\w\vp(M^j\cdot-k) \rangle$
has meaning in some sense.

For a suitable function/distribution $\phi$, the periodization of  $m^j \phi(M^j\cdot)$ leads to the sequence $\{\phi_j\}_j$  of periodic functions/distributions  such that  the $k$-th Fourier coefficient of $\phi_j$
is equal to $\h\phi({M^*}^{-j}k)$, where $\h\phi$ is the Fourier transform of $\phi$.
So, in the periodic case, the quasi-projection operators take the form

$$
Q_j(f,\phi_j,\w\phi_j)  = \frac{1}{m^j} \sum\limits_{k} \langle f, \w\phi_j(\cdot - M^{-j} k)
\rangle \phi_j (\cdot - M^{-j}k),$$
where the sum over $k$ is finite due to  the periodicity of $\phi_j$ and  the  multiplier $\frac{1}{m^j}$ is for a suitable normalization (see Section~3 for more details).
In particular, the periodization of the classical sampling expansion, where $\w\phi$ is the Dirac delta-function and $\phi$ is the sinc-function, leads to $Q_j(f,\phi_j,\w\phi_j)$, where $\widetilde \phi_j $  is a periodic distribution whose Fourier coefficients
are $\h{\widetilde \phi_j }(k)= 1$, $k\in\zd$, and  $\phi_j= m^j \Lambda_j$, where $\Lambda_j$ is a 1-periodic fundamental interpolant  on a grid $M^{-j}k$, i.e. $\Lambda_j(M^{-j}k) = \delta_{\nul k},$ $k\in  M^j [-1/2,1/2)^d \bigcap \zd.$ It is easy to see that in this case,  $Q_j(f,\phi_j,\w\phi_j)$ takes a sampling form
\begin{equation*}
  Q_j(f,\phi_j,\w\phi_j) = P_j(f)=\sum_{k} f(M^{-j}k) \Lambda_j (\cdot - M^{-j}k).
\end{equation*}

Univariate and multivariate fundamental interpolants  $\Lambda_j$  on the uniform grids and the corresponding quasi-projection operators $P_j(f)$ were investigated by many authors (see, e.g.,~\cite{Locher, Delvos, Si1, Pop, Pop2, Spr2000}).

It is well known  that in the non-periodic case  some compatibility of $\phi$ and $\w\phi$ and the Strang-Fix conditions for $\phi$ are required for successful error estimates. The periodic Strang-Fix conditions were introduced in~\cite{Brumme}. For the case of the  diagonal dilation matrix $M=2I_d$  under the periodic Strang-Fix conditions of order $s>d/2$ for the sequence $\{\Lambda_j\},$ the results in~\cite{Brumme} and~\cite{SprSick98} yield the following error estimate for $P_j$ in $L_2$ norm
$$
\|f - P_j(f)\|_{2} \le C 2^{-j \min\{\gamma,s\}} \|f\|_{H_2^{\max\{\gamma,s\}}},
$$
where $\gamma> d/2$ and $H_2^{r}$  is the Sobolev space of order $r$.

Many results in harmonic analysis involve spaces described in terms of the Fourier transform (in the non-periodic case) such as
the Fourier algebra and different  its generalizations, see, e.g.,~\cite{feich}, \cite{LST}. Similarly, in the periodic case, spaces described in terms of the Fourier coefficients appear. In particular, a natural class of such spaces, including the Sobolev space $H_2^{r}$ and the Wiener algebra, consists of the spaces ${A_q^{\alpha}} $  such that the ${A_q^{\alpha}} $-norm is a weighted $\ell_q$-norm
of the sequence of Fourier coefficients (see Section~2).
For the case $M=2I_d$,  under the strengthened Strang-Fix conditions of order $s$
the following error estimate for $P_j$  in the ${A_q^{\alpha}} $-norm  was obtained by Sprengel~\cite{Spr2000}
\be
\label{spr}
\|f - P_j(f)\|_{A_q^{\alpha}} \le C2^{-j \min\{\gamma-\alpha, s\}} \|f\|_{A_q^\gamma},
\ee
where $q \ge 1$, $\gamma \ge \alpha \ge 0$, and $\gamma > d(1-1/q)$.
For lattices generated by a matrix $M$ whose eigenvalues are greater or equal (in absolute value) than 2,
an analogous estimate was obtained by Bergmann and Prestin in~\cite{BergPrestin}.

 The goal of this paper is to obtain periodic  analogs  of some author's results in~\cite{KS, KKS, KS3}.
Namely, in~\cite{KS} the error analysis of non-periodic quasi-projection operators~(\ref{my00}) was given for a  class of tempered distributions $\w\phi$, including the Dirac delta-function, and for a wide class of functions $\phi$ with enough decay   of $\phi$ itself as well as its Fourier transform $\h\phi$. Error estimates in $L_p$-norm, $p\ge2$, were obtained  under the assumptions of the Strang-Fix conditions of order $s$ for $\phi$ and the weak compatibility of $\phi$ and $\w\phi$ of order $s$  (which means vanishing of all derivatives up to order $s$ in the origin of the function $1-\h{\w\phi}\overline{\h{\phi}}$). The obtained estimates show that the approximation order depends on the smoothness of $f$ and on $s$. In particular, it was established that the approximation order equals $s$ for smooth enough functions $f$.
Similar results were obtained in~\cite{KKS} for the same class of distributions  $\w\phi$ and a class of band-limited functions $\phi$, including the sinc-function,  under the assumption of strict compatibility of $\phi$ and $\w\phi$ (which means that the function $1-\h{\w\phi}\overline{\h{\phi}}$ is identical zero in a neighborhood of the origin). In this case, it was shown that the approximation order depends only on the smoothness of $f$.
In~\cite{KS3}, Kantorovich-Kotelnikov type operators, that are the  quasi-projection operators~(\ref{my00}) with summable $\w\phi$ and band-limited $\phi$, were investigated.   These operators are bounded in $L_p$ for $1<p<\infty$.  Under the assumption of weak compatibility of $\phi$ and $\w\phi$ of order~$s$, the $L_p$-rate of convergence was given in terms of the classical moduli of smoothness of order $s$.

The paper is organized as follows. Section~2 is devoted to notation  and basic  definitions. A wide class of periodic sampling-type quasi-projection operators $Q_j(f,\phi_j,\w\phi_j)$ with matrix dilations is introduced  in Section~3. The main results are presented in Section~4. Error estimates
in terms of the Fourier coefficients of the approximated function is given in Section~4.1.  Under the Strang-Fix conditions of order $s$ for $\phi_j$ and the weak compatibility of order $s$ of functions $\phi_j$ and distributions $\w\phi_j$, an error estimate   in the ${A_q^{\alpha}} $-norm is obtained in Theorem~\ref{theoMain}. Under the same assumptions,  an error estimate for $Q_j$  in
the $L_p$-norm,   $p\ge2$, is obtained in Theorem~\ref{big_p_th}. This theorem provides a periodic analog of the results obtained in~\cite{KS}.
In Section~4.2, under some additional assumptions on the  distributions $\w\phi_j$ and the matrix dilations $M$, we give several improvements of the error estimates obtained in the previous section. Particulary, if $\phi_j$ and  $\w\phi_j$ are strongly compatible, then we show that the error estimate for $Q_j(f,\phi_j,\w\phi_j)$ can be given only in terms of the best approximation (see Theorem~\ref{cor1}), but if $\phi_j$ and  $\w\phi_j$ are weakly compatible, then the corresponding estimates are given simultaneously in terms of Fourier coefficients and the best approximation (see Theorem~\ref{cor3}). The case of Kantorovich-type quasi-projection operators (i.e., $\w\phi_j$ is an integrable function) is also considered. In Section~5, we provide some examples.

\section{Notation}

We use the standard multi-index notations.
    Let $\n$ be the set of positive integers, $\rd$ be the $d$-dimensional Euclidean space,
    $\zd$ be the integer lattice  in $\rd$,
    $\td=\rd\slash\zd$ be the $d$-dimensional torus.
    Let  $x = (x_1,\dots, x_d)^{T}$ and
    $y =(y_1,\dots, y_d)^{T}$ be column vectors in $\rd$,
    then $(x, y):=x_1y_1+\dots+x_dy_d$,
    $|x| := \sqrt {(x, x)}$; $\nul=(0,\dots, 0)^T\in \rd$;  		
		$\z_+^d:=\{x\in\zd:~x_k\geq~{0}, k=1,\dots, d\}.$
		If $a\in\rd$, $r>0$, then $B_r(a)$
		denotes the ball of radius $r$ with the center in $a$.

    If $\alpha\in\zd_+$,  we set
    $[\alpha]=\sum\limits_{k=1}^d \alpha_k$,
$D^{\alpha}f=\frac{\partial^{[\alpha]} f}{\partial^{\alpha_1}x_1\dots
	\partial^{\alpha_d}x_d}$.

If $A$ is a $d\times d$ matrix,
then $\|A\|$ denotes its operator norm in $\rd$; $A^*$ denotes the conjugate matrix to $A$.

Let $M$ be a dilation matrix, i.e. an integer valued $d\times d$ matrix, such that the absolute value of each its eigenvalue is greater than 1, $m:=|\det M|$,   $D(M) : = M [-1/2,1/2)^d \cap \zd$.
It is known (see, e.g.,~\cite[Chapter 2]{KPS}) that  $D(M)$ is a set of digits
of $M$, and any $k\in\zd$ can be uniquely represented as $k=M n + r$, $r\in D(M)$, $n\in\zd.$

Since the spectrum of the operator $M^{-1}$ is
located in $B_r(\nul)$, where $r=r(M^{-1}):=\lim\limits_{j\to+\infty}\|M^{-j}\|^{\frac1j}$ is
the spectral radius of $M^{-1}$, and there exists at least
one point of the spectrum on the boundary of the ball, we have
	\be
	\|M^{-j}\|\le {C_{M,\vartheta}}\, \vartheta^{-j},\quad j\ge0,
	\label{00}
	\ee
for every  positive number $\vartheta$  whose absolute value is smaller
than absolute value of any eigenvalue of  $M$.
In particular, we can take $\vartheta > 1$ and, hence,
	$	\lim_{j\to+\infty}\|M^{-j}\|=0.	$
	
A  matrix $M$ is called isotropic
	if it is similar to a diagonal matrix
such that numbers $\lambda_1,\dots,\lambda_d$ are placed on the main diagonal
	and $|\lambda_1|=\cdots=|\lambda_d|$.
	Thus, $\lambda_1,\dots,\lambda_d$ are eigenvalues of $M$
	and the spectral radius of $M$ is equal to $|\lambda|,$
	where $\lambda$ is one of the eigenvalues of $M.$
	Note that if the matrix $M$ is isotropic then
	$M^*$ and $M^j$ are isotropic for all $j\in\z.$	
	It is well known that for an isotropic matrix $M$ and for  any $j\in\z$
	we have
	\be
	C^M_1 |\lambda|^j \le \|M^j\| \le C^M_2 |\lambda|^j,\quad j\in\z,
	\label{10}
\ee
   where $\lambda$ is one of the eigenvalues of $M,$ constants $C^M_1, C^M_2$ do not depend on $j$.

 We  will use notation $L_p$ for the space  $L_p(\td)$ with the usual norm
$\|f\|_p = \Big(\int_{\td}|f(x)|^p dx\Big)^{1/p}$ for $1\le p<\infty$, and $\|f\|_\infty={\rm vrai}\sup |f|$.

If $f\in L_1$, then $\h f(k)$, $k\in\zd$, denotes the $k$-th Fourier coefficient of $f$.

 For $g\in L_1$, we will use the following  notation $g^-(t)=\overline{g(-t)}$.

 As usual, the convolution of appropriate functions $f$ and $g$ is given by
$$(f*g)(x)=\int_{\T^d} f(x-t)g(t)dt.$$

We  will use notation $\theta_\alpha(x): = (1+ |x|^2)^{\frac{\alpha}{2}}$, $\alpha \ge 0,$ $x\in \rd.$
It is easy to check that
\be
\label{103}
\theta_\alpha(x+y) \le \theta_\alpha(x) \theta_\alpha(y)\quad\forall x,y\in\rd.
\ee
Since $\|M\|> 1$  for any  dilation matrix $M$, and hence $\|M^{*j}\|> 1$  for all $j\in\n$, we have
\begin{equation}\label{103++}
  \theta_\alpha(x) = \theta_\alpha(M^{*j} M^{*-j}x) \le \|M^{*j}\|^{\alpha} \theta_\alpha(M^{*-j}x).
\end{equation}
Also, for $\alpha \ge 0$ and $1 \le q \le \infty$, we will use notation
$$
\ell_q^\alpha : = \left\{ x = \{x_k\}_k: \{\theta_\alpha(k) x_k\}_k \in \ell_q\right\},\quad \ell_q=\ell_q(\zd),
$$
and  $\|x\|_{\ell_q^\alpha} = \| \{\theta_\alpha(k) x_k\}_k\|_{\ell_q}$ for $x\in \ell_q^\alpha$.

Let $1 \le q \le \infty$  and $\alpha \ge0$. The space $A_q^\alpha$ of periodic functions is defined by
$$
A_q^\alpha: = \left\{ f\in L_1: \{\h f(k)\}_k\in \ell_q^\alpha
\right\}, \quad \|f\|_{A_q^\alpha}:=\|\{\h f(k)\}_k\|_{\ell_q^\alpha}.
$$
These spaces include the Wiener algebra $A:=A_1^0$ of functions with absolutely convergent Fourier series.
It is clear that $A_q^\beta \subset A_q^\alpha$ for $\beta \ge \alpha.$ If $q=2$, then the space $A_2^\alpha$ coincides with the Sobolev space $H_2^{\alpha}.$ So $\alpha$ can be considered as a smoothness parameter. However, these smoothness
properties (except $p=2$) differ from usual smoothness of the fractional Sobolev spaces $H_p^{\alpha}$. Some embeddings of $A_q^\alpha$ into $H_p^\beta$ and vice versa can be found in~\cite{Spr2000}.

For a function $f\in A_q^\alpha$, we set
$$ \|f\|^{Out}_{A_q^\alpha,j}:=\Big(\sum_{k\notin D(M^{*j})} |\h f(k)|^q
\theta^q_{\alpha}(k)\Big)^{1/q},
\quad \|f\|^{In}_{A_q^\alpha,j}:=\Big(\sum_{k\in D(M^{*j})} |\h f(k)|^q
\theta^q_{\alpha}(k) \Big)^{1/q}.
$$
for  $1\le q<\infty$, and
		$$ \|f\|^{Out}_{A_\infty^\alpha,j}:=\sup\limits_{k\notin D(M^{*j})} |\h f(k)|
	\theta_{\alpha}(k),
	\quad \|f\|^{In}_{A_\infty^\alpha,j}:=\sup\limits_{k\in D(M^{*j})} |\h f(k)|
	\theta_{\alpha}(k) 
	$$
for $q=\infty$.
For  convenience, we will also use the following notation:
$$ \|f\|^{q,Out}_{A_q^\alpha,j}:=\sum_{k\notin D(M^{*j})} |\h f(k)|^q
\theta^q_{\alpha}(k),
\quad \|f\|^{q,In}_{A_q^\alpha,j}:=\sum_{k\in D(M^{*j})} |\h f(k)|^q
\theta^q_{\alpha}(k). 
$$

Let ${\cal D} = C^{\infty}(\td)$ be the space of infinitely differentiable functions on $\rd$ that are periodic with period 1. A continuous linear functional on the space ${\cal D}$ is a periodic distribution.
The linear space of periodic distributions we denote by ${\cal D}'$. For a periodic distribution $\phi\in {\cal D}'$ and a function $f\in {\cal D}$, we denote the action of $\phi$ on $f$ by $\phi(f)$. For convenience, we will use notation $\langle \phi, f \rangle : = \phi(\overline{f})$ and $\langle f, \phi \rangle = \overline{\langle \phi, f \rangle}$.
It is known (see, e.g.,~\cite[p. 322]{Folland},~\cite[p. 144]{Triebel}) that any periodic distribution $\phi$ can be expanded in a weakly convergent (in ${\cal D}'$) Fourier series
\begin{equation}
\phi(x) = \sum_{n\in\zd} \h\phi(n) e^{2\pi i (n, x)},
\label{fPDisrt}
\end{equation}
where the sequence $\{\h\phi(n)\}_n$ has at most polynomial growth and for any $f \in {\cal D}$
$$
\langle  f, \phi \rangle = \sum_{n\in\zd} \h f(n) \overline{\h\phi(n)}.
$$
Also, conversely, for any sequence $\{\h\phi(n)\}_n$ of at most polynomial growth the series in the right-hand side of~(\ref{fPDisrt}) converges weakly to a periodic distribution.
The numbers $\h\phi(n)$ are called the Fourier coefficients of a periodic distribution $\phi$ and
$\h\phi(n) = \langle \phi, e^{2\pi i (n, \cdot)} \rangle = \phi(e^{-2\pi i (n, \cdot)})$. The convolution of $f\in \mathcal{D}$ and the distribution $\vp$ is defined by
$
(f*\vp) (x)=\langle f, \overline{\vp(x-\cdot)}\rangle.
$


\section{Sampling-type quasi-projection operators}

In this section, we define  the periodic quasi-projection operators $Q_j(f,\phi_j,\w\phi_j)$, where $\phi\in L_1$, $\w\phi_j\in {\cal D}'$, $j\in\n$, by
$$
Q_j(f,\phi_j,\w\phi_j)  = \frac{1}{m^j} \sum\limits_{k\in D(M^j)} \langle f, \w\phi_j(\cdot - M^{-j} k)
\rangle \phi_j (\cdot - M^{-j}k).
$$
By the definition of ${\cal D}'$, such an operator $Q_j$ has meaning only for
 $f \in {\cal D}$, and
\be
\label{100my}
 \langle f, \w\phi_j(\cdot - M^{-j} k) \rangle
 =  \sum_{l\in\zd}
 \h f(l) \overline{\h{\w\phi_j}(l)}   e^{2 \pi i (k,M^{*-j}l)}.
\ee
But since  the Fourier coefficients of $\w\phi_j$ have polynomial growth, the latter series converges  for any $f$ whose Fourier coefficients decay sufficiently fast.
Thus, to extend the class of functions $f$, we define
the inner product  $ \langle f, \w\phi_j(\cdot - M^{-j} k) \rangle$  by~(\ref{100my}).

Next, we introduce several conditions on the sequences of functions $\phi_j$
and distributions $\w \phi_j$
under which the quasi-projection operator $Q_j(f,\phi_j,\w\phi_j)$ provides good enough approximation of an appropriate function $f$.
For additional motivation of these conditions, we  consider their connection with the analogous conditions for the non-periodic case (see~\cite{KS, KKS}). For this, we note that the periodic and non-periodic cases can be connected via periodization as follows. Let  $\phi \in L_1({\mathbb R}^d)$ with $\{\widehat\phi(k)\}_{k\in\zd} \in \ell_1$. Setting $\phi_{j,l}(x):=m^{j/2}\phi(M^jx+l)$  and
\begin{equation}
\phi_j(x) := m^{j} \sum_k  \phi(M^{j}(x+k)) = m^{j/2} \sum_k \phi_{j,0}(x+k),
\label{fPeriodiz}
\end{equation}
where the multiplier $ m^{j/2}$ is for a suitable normalization, we have that each $\phi_j$ is in $L_1$ and  by the Poisson summation formula,
$$\phi_j(x) =  m^{j/2} \sum_k \h{\phi_{j,0}}(k) e^{2\pi i (k,x)} =
\sum_k \h \phi(M^{*-j}k)e^{2\pi i (k,x)} = \sum_k   \h{\phi_j}(k) e^{2\pi i (k,x)},
$$
which implies that  $\h{\phi_j}(k) = \h \phi(M^{*-j}k).$ 

\smallskip

    1. The conditions on the growth of order $N\ge 0$ for Fourier coefficients of $\{\w\phi_j\}_j$:
\begin{equation}
\label{fWPhi}
 \begin{split}
| \h{\w\phi_j}(k)| & \le C_{\w\phi} |M^{*-j}k|^{N}, \quad  \forall k \notin D(M^{*j}), \quad \forall j\in\n,
\\
\max\limits_{k\in D(M^{*j})} |  \h{\w\phi_j}(k)| & \le  C_{\w\phi}, \phantom{|M^{*-j}k|^{N}} \quad \forall j\in\n.
 \end{split}
\end{equation}
These conditions correspond to the following non-periodic conditions:
$$
|\h{\w\phi}(\xi)|\le C_{\w\phi}|\xi|^N\quad \text{for}\quad \xi \notin [-1/2,1/2]^d
$$
and
$$
|\h{\w\phi}(\xi)|\le C_{\w\phi}\quad \text{for}\quad \xi \in [-1/2,1/2]^d.
$$
To show this, we note that 
$ \h{\w\phi_j}(k) = \h{\w\phi}(M^{*-j} k)$. Thus, if non-periodic conditions are valid for $\w\phi$, then conditions~(\ref{fWPhi}) are valid for $\w\phi_j$.

2. The Strang-Fix conditions of order $s$, $s>0$, for a sequence of functions $\{\phi_j\}_j$:
\begin{equation}
| \h{\phi_j}(M^{*j}n+r)| \le b_n |M^{*-j} r|^s, \quad \forall n\neq\nul, \quad \forall r\in D(M^{*j}).
\label{fPhiSF}
\end{equation}
This  corresponds to the well known  Strang-Fix conditions for a non-periodic function $\phi$: $D^{\beta}\h{\phi}(n) = 0$, for $n\neq \nul,$ $[\beta]<s$, $s\in\n$ and the condition that $\h{\phi}$ is boundedly differentiable up to order $s$. Indeed, by Taylor's formula near the point $\xi = n$, we have
for $r\in D(M^{*j})$
$$
\h{\phi}(n + M^{*-j} r) = \sum_{[\beta]=s} \frac{s}{\beta!}  (M^{*-j} r)^{\beta} \int_0^1 (1-t)^{s-1} D^{\beta}\h{\phi}(n + t M^{*-j} r) d t.
$$
Thus, since $  \h{\phi_j}(M^{*j}n+r)= \h{\phi}(n + M^{*-j} r)$, we get 
$| \h{\phi_j}(M^{*j}n+r)| \le C_{s,\phi} |M^{*-j} r|^s.$

For  error estimates in the non-periodic case (see, e.g.,~\cite{KS}), the following additional assumption helps:
$$
\sum_{n\in\zd}|D^{\beta}\h{\phi}(n +\xi)|^q\le B_s,\quad \forall\xi\in\rd, [\beta]<s.
$$
A periodic analog  of this condition is $\{b_n\}_n\in \ell_q$.

3. The  weak compatibility conditions  of order $s$, $s> 0$, for $\{\phi_j\}_j$ and $\{\w \phi_j\}_j$:
\begin{equation}
\left|1 -  \h{\phi_j}(r)\overline{\h{\w\phi_j}(r)}  \right| \le b_0  |M^{*-j} r|^s, \quad \forall r\in D(M^{*j}), \quad \forall j\in\n.
\label{fPhiPhiW}
\end{equation}
This  corresponds to the following non-periodic condition:  $D^{\beta}(1-\h{\phi}\overline{\h{\w\phi}})(\nul) = 0$, $[\beta]<s$,
$s\in\n$, and condition that $\h{\phi}\overline{\h{\w\phi}}$ is boundedly differentiable up to order $s$. Indeed, by Taylor's formula near the point $\xi = 0$,
\begin{equation}\label{11+}
  \h{\phi} (M^{*-j} r)\overline{\h{\w\phi}(M^{*-j} r)} = 1 + \sum_{[\beta]=s} \frac{s}{\beta!}  (M^{*-j} r)^{\beta} \int_0^1 (1-t)^{s-1} D^{\beta}\h{\phi}\overline{\h{\w\phi}}( t M^{*-j} r) d t.
\end{equation}
Thus, since $\h{\phi}(M^{*-j} r) \overline{\h{\w\phi}(M^{*-j} r)} = \h{\phi_j}(r)\overline{\h{\w\phi_j}(r)} $, we have
$
| 1 -   \h{\phi_j}(r)\overline{\h{\w\phi_j}(r)} | \le C_{s,\w\phi} |M^{*-j} r|^s.
$

4. The conditions on the uniform boundedness of the Fourier coefficients of $\{\phi_j\}_j$:
\begin{equation}
\max\limits_{k\in D(M^{*j})} |\h{\phi_j}(k)| \le C_{\phi}, \quad \forall j\in\n.
\label{fPhi}
\end{equation}

\section{Error estimates for $Q_j(f,\vp_j,\w\vp_j)$}
\subsection{Estimates in terms of the Fourier coefficients}

First we establish two utility lemmas.
For convenience, in the proofs of these lemmas we will use the following notation
$$
\label{fCd}
C_d :=  \max\limits_{x\in[-\frac 12,\frac 12]^d} \theta_1(x) = \left(1+\frac d4\right)^{\frac 12}.
$$

\begin{lem}
	\label{lem1}
	Let $1\le q \le \infty$, $1/p+1/q=1$, $j\in\n$, $\alpha\ge 0$,
	$\{\w\phi_j\}_j$ be a sequence of periodic distributions
	whose Fourier coefficients satisfy condition~(\ref{fWPhi}) with parameter $N\ge 0$.
	Suppose $f\in A_q^{\gamma}$, where $\gamma>N+ d/p$ for $q\ne1$ or  $\gamma \ge N$ for $q=1$.
	Then, if $1 \le q <\infty$
	\begin{equation*}
	\left(\sum\limits_{r\in\,D(M^{*j})}
	\left| \sum_{n\neq 0} \h f(r + M^{*j}n)  \overline{\h{\w\phi_j}(r + M^{*j}n)}
	\right|^q \theta^q_\alpha(r) \right)^{\frac 1q} \le
	C_{Lem} \|{M^{*j}}\|^{\alpha} \|M^{*-j}\|^{ \gamma}  
	\|f\|^{Out}_{A_q^{\gamma},j},
	\end{equation*}
	and if $q=\infty$
	\begin{equation*}
	\sup\limits_{r\in\,D(M^{*j})}
	\left| \sum_{n\neq 0} \h f(r + M^{*j}n)  \overline{\h{\w\phi_j}(r + M^{*j}n)}
	\right| \theta_\alpha(r) \le
	C_{Lem} \|{M^{*j}}\|^{\alpha} \|M^{*-j}\|^{ \gamma}  
	\|f\|^{Out}_{A_\infty^{\gamma},j},
	\end{equation*}
	where the constant $C_{Lem}$ does not depend on $j$ and $f$. 
\end{lem}

{\bf Proof.} First, let $1<q<\infty$ and $\alpha=0$.
Applying H\"older's inequality, condition~\eqref{fWPhi}, and the fact that $\frac{1}{|M^{*j}k|} \le \frac{\|M^{*-j}\|}{|k|}$, we derive
\begin{equation*}
\begin{split}
&\sum\limits_{r\in\,D(M^{*j})}
\left| \sum_{n\neq \nul} \h f(r + M^{*j}n) \overline{\h{\w\phi_j}(r + M^{*j}n)}
\right|^q \\
& \hspace{1cm}\le \sum\limits_{r\in\,D(M^{*j})}
\left( \sum_{n\neq \nul} \frac{1}{|r + M^{*j}n|^{p (\gamma-N)}}
\right)^{\frac qp}\\
&\hspace{1cm}\hspace{1cm}\hspace{1cm}\hspace{1cm}\hspace{1cm}\hspace{1cm}\times\sum_{n\neq \nul} |\h f(r + M^{*j}n)|^q |\h{\w\phi_j}(r + M^{*j}n)|^q |r + M^{*j}n|^{q (\gamma-N)} \\
&\hspace{1cm} \le \max\limits_{r\in\,D(M^{*j})} \left( \sum_{n\neq \nul} \frac{1}{|r + M^{*j}n|^{p (\gamma-N)}}
\right)^{\frac qp}
\sum\limits_{k\notin D(M^{*j})} |\h f(k)|^q |\h{\w\phi_j}(k)|^q |k|^{q (\gamma-N)} \\
&\hspace{1cm}\le \|M^{*-j}\|^{(\gamma-N)q} \max\limits_{\xi\in\,[-1/2,1/2)^d} \left( \sum_{n\neq \nul} \frac{1}{|n+\xi|^{p (\gamma-N)}}
\right)^{\frac qp}  C^q_{\w\phi} \|M^{*-j}\|^{qN}	\sum\limits_{k\notin D(M^{*j})} |\h f(k)|^q    |k|^{q \gamma} \\
&\hspace{1cm} \le C^q_{p,\gamma,N} C^q_{\w\phi} \|M^{*-j}\|^{q \gamma}  
\|f\|^{q,Out}_{A_q^\gamma,j}.
\end{split}
\end{equation*}
Repeating the same steps, the  required  estimate for $q = \infty$ can be derived.
For $q=1$, we have
\begin{equation*}
\begin{split}
&\sum\limits_{r\in\,D(M^{*j})}
\left| \sum_{n\neq \nul} \h f(r + M^{*j}n) \overline{\h{\w\phi_j}(r + M^{*j}n)}
\right| \le \sum\limits_{k\notin\,D(M^{*j})}
|\h f(k) \overline{\h{\w\phi_j}(k)}| \\
&\hspace{1cm}  \le  C_{\w\phi} \sum\limits_{k\notin\,D(M^{*j})}
|\h f(k)|  |M^{*-j}k|^{N}  \le
2^\gamma C_{\w\phi}  \sum\limits_{k\notin\,D(M^{*j})}
|\h f(k)| |M^{*-j}k|^{\gamma}\le
2^\gamma C_{\w\phi} \|M^{*-j}\|^\gamma  
\|f\|^{Out}_{A_1^\gamma,j}
\end{split}
\end{equation*}
since $|M^{*-j}k| \ge \frac 12$ for $k\notin\,D(M^{*j})$ and $|M^{*-j}k|^{N} \le (2 |M^{*-j}k|)^N \le (2 |M^{*-j}k|)^\gamma.$

To prove the lemma for $\alpha>0$, it is sufficient to note that
\begin{equation}
\max_{r\in D(M^{*j})} \theta^q_\alpha(r) \le \|M^{*j}\|^{q\alpha} \max_{r\in D(M^{*j})}\theta^q_\alpha(M^{*-j}r)  \le \|M^{*j}\|^{q\alpha}  \max\limits_{x\in[-\frac 12,\frac 12]^d} \theta_\alpha^q(x)  \le C_d^{\alpha q}
\|M^{*j}\|^{q\alpha},
\label{fLemUtil}
\end{equation}
and $C_{Lem} = C_{p,\gamma,N} C_{\w\phi} C^\alpha_d$ for $q>1$ and $C_{Lem} = 2^\gamma C_{\w\phi} C^\alpha_d$ for $q=1.$ $\Diamond$

\begin{rem}
The constant $C_{p,\gamma,N}$ in the proof of  Lemma~\ref{lem1} can be estimated  as follows (see \cite[Lemma~1.10]{Pop})
$$
C^p_{p,\gamma,N} = \max\limits_{\xi\in\,[-\frac 12,\frac 12)^d}  \sum_{n\neq \nul} \frac{1}{|n+\xi|^{p (\gamma-N)}} \le  2^{p (\gamma-N)} \sum_{v=1}^d 2^v \frac{1}{v^{p (\gamma-N) / 2}} \binom{d}{v} \left( 1 + \frac{v}{2(p (\gamma-N) - v)} \right)^v.
$$
\end{rem}

\begin{lem}
\label{lem2}
Let $1\le q \le\infty$, $1/p+1/q=1$, $j\in\n$, $\alpha\ge 0$,
$\{\phi_j\}_j$ and $\{\w\phi_j\}_j$ be sequences of periodic functions and periodic distributions, respectively,
whose Fourier coefficients are such that conditions (\ref{fWPhi}), (\ref{fPhiSF}), (\ref{fPhiPhiW}), (\ref{fPhi}) are valid with parameters $N\ge 0$, $s>0$, and $b=\{b_n\}_n \in \ell_q^\alpha$.
 Suppose $f\in A_q^{\gamma}$, where $\gamma\ge \alpha$ and $\gamma>N+ d/p$ for $q\ne1$ or  $\gamma\ge N$ for $q=1$,
$$
I:=
\begin{cases}
\displaystyle\Big(\sum\limits_{n\in\zd} \Big| \h f(n) - \sum_{l\in\zd} \h f(n + M^{*j}l)
	\overline{\h{\w\phi_j}(n + M^{*j}l)}  \h{\phi_j}(n) \Big|^q 
	\theta^q_\alpha(n)  \Big)^{\frac1q},\  &{\rm{if}}\ q<\infty,
\\
\displaystyle\sup\limits_{n\in\zd} \Big| \h f(n) - \sum_{l\in\zd} \h f(n + M^{*j}l)
	\overline{\h{\w\phi_j}(n + M^{*j}l)}  \h {\phi_j}(n) \Big|
	\theta_\alpha(n),\  &{\rm{if}}\  q=\infty.
\end{cases}
$$
 Then
\begin{equation}
I \le
	C_1     \|{M^{*j}}\|^{\alpha} \|M^{*-j}\|^{ s}\|f\|^{In}_{A_q^{s+\alpha},j} + C_2 \  \|{M^{*j}}\|^{\alpha} \|M^{*-j}\|^{\gamma} \|f\|^{Out}_{A_q^\gamma,j},
\label{lemMain}
\end{equation}
where the constants $C_1$ and $C_2$ do not depend on $j$ and $f$. 

\end{lem}

{\bf Proof.} First, we consider the case $q<\infty$. Set
$$
J:=I^q=\sum_{n\in\zd} \left| \h f(n) - \sum_{l\in\zd} \h f(n + M^{*j}l)
	\overline{\h{\w\phi_j}(n + M^{*j}l)}  \h{\phi_j}(n) \right|^q
	\theta^q_\alpha(n)
$$
and split the sum $J$ into two parts such that
$$
J=\sum_{n\in D(M^{*j})}\, + \sum_{n \notin D(M^{*j})}:=J_0+J_1.
$$
%
%
Estimating $J_0$, we derive
\begin{equation*}
\begin{split}
 J_0
	 = &
     \sum_{n\in D(M^{*j})} \left| \h f(n) - \h f(n)
	\overline{\h{\w\phi_j}(n)}  \h{\phi_j}(n) -
    \sum_{l\neq \nul} \h f(n + M^{*j}l)
	\overline{\h{\w\phi_j}(n + M^{*j}l)}  \h{\phi_j}(n) \right|^q \theta^q_\alpha(n)  \\
	\le &
	 2^{\frac qp} \sum_{n\in D(M^{*j})} |\h f(n) - \h f(n)
	\overline{\h{\w\phi_j}(n)}  \h{\phi_j}(n)|^q \theta^q_\alpha(n)
	\\
	 & \hspace{1cm} + 2^{\frac qp} \sum_{n\in D(M^{*j})} \left|\sum_{l\neq \nul} \h f(n + M^{*j}l)
	\overline{\h{\w\phi_j}(n + M^{*j}l)}  \h{\phi_j}(n) \right|^q \theta^q_\alpha(n) =:	J_{00} + J_{01}.
\end{split}
\end{equation*}
For the first term $J_{00}$, by the compatibility conditions for $\phi_j$ and $\w \phi_j$ of order $s$, see~\eqref{fPhiPhiW}, we get
\begin{equation*}
\begin{split}
J_{00} = & 2^{\frac qp} \sum_{n\in D(M^{*j})} |\h f(n)|^q | 1 -
	\overline{\h{\w\phi_j}(n)}  \h{\phi_j}(n)|^q \theta^q_\alpha(n) \\
	\le &
2^{\frac qp} b_0^q \|M^{*-j}\|^{q s}  \sum_{n\in D(M^{*j})}  |\h f(n)|^q |n|^{qs} \theta^q_\alpha(n) \le
2^{\frac qp} b_0^q \|M^{*-j}\|^{q s} \|f\|^{q,In}_{A_q^{s+\alpha},j}.
\end{split}
\end{equation*}
For the second term  $J_{01}$, by Lemma~\ref{lem1} and~\eqref{fPhi}, we have
\begin{equation*}
\begin{split}
J_{01} & = 2^{\frac qp}  \sum_{n\in D(M^{*j})} \left|\sum_{l\neq \nul} \h f(n + M^{*j}l)
	\overline{\h{\w\phi_j}(n + M^{*j}l)}  \h{\phi_j}(n) \right|^q \theta^q_\alpha(n) \\
	& \le
2^{\frac qp} \max\limits_{n\in D(M^{*j})} |\h{\phi_j}(n)|^q \sum_{n\in D(M^{*j})} \left|\sum_{l\neq \nul} \h f(n + M^{*j}l)
	\overline{\h{\w\phi_j}(n + M^{*j}l)}\right|^q \theta^q_\alpha(n) \\
	& \le
2^{\frac qp}  C^q_{\phi} C^q_{Lem} \|{M^{*j}}\|^{q\alpha} \|M^{*-j}\|^{q\gamma} \|f\|^{q,Out}_{A_q^{\gamma},j}.
\end{split}
\end{equation*}
Thus, the term $J_0$ is estimated by
\be
J_0 \le 2^{\frac qp} \left(b_0^q  \|M^{*-j}\|^{q s} \|f\|^{q,In}_{A_q^{s+\alpha},j} +  C^q_{\phi} C^q_{Lem} \|{M^{*j}}\|^{q\alpha} \|M^{*-j}\|^{q\gamma} \|f\|^{q,Out}_{A_q^{\gamma},j}\right).
\label{J0}
\ee

Consider $J_1.$ By Minkowski's inequality, we get
\begin{equation*}
\begin{split}
	 J_1 & : =  \sum_{n\notin D(M^{*j})}  \left| \h f(n) - \h f(n)
	\overline{\h{\w\phi_j}(n)}  \h{\phi_j}(n) -
    \sum_{l\neq \nul} \h f(n + M^{*j}l)
	\overline{\h{\w\phi_j}(n + M^{*j}l)} \h{\phi_j}(n) \right|^q  \theta^q_\alpha(n) \\
	& \le
	 2^{\frac qp} \sum_{n\notin D(M^{*j})}   \left| \h f(n) - \h f(n)
	\overline{\h{\w\phi_j}(n)}  \h{\phi_j}(n) \right|^q \theta^q_\alpha(n)
	\\
	& \hspace{1cm} + 2^{\frac qp} \sum_{n\notin D(M^{*j})} \left|
    \sum_{l\neq \nul} \h f(n + M^{*j}l)
	\overline{\h{\w\phi_j}(n + M^{*j}l)}  \h{\phi_j}(n) \right|^q \theta^q_\alpha(n)   =: J_{10} + J_{11}.
\end{split}
\end{equation*}
Again, by Minkowski's inequality,
	$$
    J_{10} \le 2^{\frac{2q}{p}}	\sum_{n\notin D(M^{*j})} |\h f(n)|^q \theta^q_\alpha(n) + 2^{\frac{2q}{p}}\sum_{n\notin D(M^{*j})}
    \left| \h f(n)
	\overline{\h{\w\phi_j}(n)}  \h{\phi_j}(n) \right|^q \theta^q_\alpha(n) = : J_{100} + J_{101}.
$$
For the sum $J_{100}$, using the inequality $\frac{1}{|n|} \le \frac{\|M^{*-j}\|}{|M^{*-j} n|}$, we derive
\begin{equation*}
\begin{split}
J_{100} &= 2^{\frac{2q}{p}}	\sum_{n\notin D(M^{*j})} \frac{|n|^{q(\gamma-\alpha)}}{|n|^{q(\gamma-\alpha)}}|\h f(n)|^q \theta^q_\alpha(n)  \le   2^{\frac{2q}{p}}\|f\|^{q,Out}_{A_q^{\gamma},j} \max\limits_{n \notin D(M^{*j})} \frac{1}{|n|^{q(\gamma-\alpha)}} \\
 &\le
   2^{\frac{2q}{p}}\|f\|^{q,Out}_{A_q^{\gamma},j} \|M^{*-j}\|^{q(\gamma-\alpha)} \max \limits_{n \notin D(M^{*j})} \frac{1}{|M^{*-j} n|^{q(\gamma-\alpha)}}\\
 &\le   2^{\frac{2q}{p}} 2^{q (\gamma-\alpha)} \|M^{*-j}\|^{q(\gamma-\alpha)} \|f\|^{q,Out}_{A_q^{\gamma},j},
\end{split}
\end{equation*}
where the last inequality is valid since for $n \notin D(M^{*j})$ we have $|M^{*-j} n| > 1/2$.

Consider $J_{101}$. By the Strang-Fix conditions~(\ref{fPhiSF}) and inequalities~(\ref{103}), (\ref{103++}) and~(\ref{fLemUtil}), for any  $n=M^{*j}k+r$ with $k\neq \nul$ and $r\in D(M^{*j})$,  we obtain
$$
| \h{\phi_j}(n)|^q  \theta^q_\alpha(n) \le b^q_k |M^{*-j} r|^{qs} \|M^{*j}\|^{\alpha q}  \theta^q_\alpha(k) \theta^q_\alpha(M^{*-j}r)
\le \|b\|_{\ell^\alpha_\infty}^q   C_d^{q(s+\alpha)}
\|M^{*j}\|^{\alpha q}.
$$
Therefore, 
\begin{equation*}
\begin{split}
J_{101} & = 2^{\frac{2q}{p}} \sum_{n\notin D(M^{*j})}
    \left| \h f(n)
	\overline{\h{\w\phi_j}(n)}  \h{\phi_j}(n) \right|^q  \theta^q_\alpha(n) \\
& = 2^{\frac{2q}{p}} \|b\|_{\ell^\alpha_\infty}^q   C_d^{q(s+\alpha)}     \|M^{*j}\|^{\alpha q}
  \sum_{n\notin D(M^{*j})}
    \left| \h f(n)
	\overline{\h{\w\phi_j}(n)} \right|^q  \\
	& \le
	  2^{\frac{2q}{p}}  \|b\|_{\ell^\alpha_\infty}^q   C_d^{q(s+\alpha)}      \|M^{*j}\|^{\alpha q} C^q_{\w\phi} \|M^{*-j}\|^{Nq} \sum_{n\notin D(M^{*j})}
    | \h f(n)  |^q  |n|^{q N} \\
&  = 2^{\frac{2q}{p}} \|b\|_{\ell^\alpha_\infty}^q   C_d^{q(s+\alpha)}     C^q_{\w\phi} \|M^{*j}\|^{\alpha q}   \|M^{*-j}\|^{Nq}\sum_{n\notin D(M^{*j})}
    | \h f(n)  |^q   \frac{|n|^{q \gamma}}{ |n|^{q (\gamma-N)}}\\
  & \le 2^{\frac{2q}{p}} \|b\|_{\ell^\alpha_\infty}^q   C_d^{q(s+\alpha)}    C^q_{\w\phi}   \|M^{*j}\|^{\alpha q}   \|M^{*-j}\|^{Nq} \sum_{n\notin D(M^{*j})}
    | \h f(n)  |^q  |n|^{q \gamma} \frac{\|M^{*-j}\|^{q (\gamma-N)}}{ |M^{*-j}n|^{q (\gamma-N)}} \\
    & \le
    2^{\frac{2q}{p}} 2^{q (\gamma-N)} \|b\|_{\ell^\alpha_\infty}^q   C_d^{q(s+\alpha)}      C^q_{\w\phi}
    \|M^{*j}\|^{\alpha q}
    \|M^{*-j}\|^{q\gamma} \|f\|^{q,Out}_{A_q^{\gamma},j},
\end{split}
\end{equation*}
where the last inequality is valid because $|M^{*-j} n| > 1/2$  whenever $n\notin D(M^{*j})$, and, therefore, 
$$
\frac{1}{ |M^{*-j}n|^{q (\gamma-N)} }\le \max\limits_{n\notin D(M^{*j})} \frac{1}{ |M^{*-j}n|^{q (\gamma-N)}}= 2^{q (\gamma-N)}.
$$

Combining the estimates for $J_{100}$ and $J_{101}$, we get
$$
J_{10} \le   2^{\frac{2q}{p}} \left(  2^{q (\gamma-\alpha)}  + 2^{q (\gamma-N)}  \|b\|_{\ell^\alpha_\infty}^q   C_d^{q(s+\alpha)}    C^q_{\w\phi}\right) \|{M^{*j}}\|^{q\alpha} \|M^{*-j}\|^{q\gamma} \|f\|^{q,Out}_{A_q^{\gamma},j}.
$$

Now, let us estimate $J_{11}$. By Minkowski's inequality,
\begin{equation*}
\begin{split}
J_{11} &=2^{\frac{q}{p}}\!\!\!\!   \sum_{n\notin D(M^{*j})}   \left|
    \sum_{l\neq \nul} \h f(n + M^{*j}l)
	\overline{\h{\w\phi_j}(n + M^{*j}l)}  \h{\phi_j}(n) \right|^q \theta^q_\alpha(n) \\
	&=   2^{\frac{q}{p}}\!\!  \sum_{r\in D(M^{*j})} \sum_{k\neq\nul}  |\h{\phi_j}(r + M^{*j}k) |^q \left|
    \sum_{l\neq \nul} \h f(r + M^{*j}(k+l))
	\overline{\h{\w\phi_j}(r + M^{*j}(k+l))} \right|^q
	\theta^q_\alpha(r + M^{*j}k) \\
  & =  2^{\frac{2q}{p}} \!\!\!\!  \sum_{r\in D(M^{*j})} \sum_{k\neq\nul}  | \h{\phi_j}(r + M^{*j}k) |^q
  |\h f(r )
	\overline{\h{\w\phi_j}(r )}|^q \theta^q_\alpha(r + M^{*j}k) \\
	& +
	2^{\frac{2q}{p}}\!\!\!\! \sum_{r\in D(M^{*j})} \sum_{k\neq\nul}  | \h{\phi_j}(r + M^{*j}k) |^q
    \left| \sum_{\underset{l\neq -k}{l\neq\nul}} \h f(r + M^{*j}(k+l))
	\overline{\h{\w\phi_j}(r + M^{*j}(k+l))} \right|^q \!\!\theta^q_\alpha(r + M^{*j}k)\\
&=:J_{110} + J_{111}
\end{split}
\end{equation*}
 Using the Strang-Fix conditions~\eqref{fPhiSF}, inequalities~\eqref{103}, \eqref{103++} and conditions~\eqref{fWPhi}, we obtain
\begin{equation*}
\begin{split}
J_{110}
	& \le
		2^{\frac{2q}{p}}  \sum_{r\in D(M^{*j})}  |M^{*-j} r|^{sq}
  |\h f(r )
	\overline{\h{\w\phi_j}(r )}|^q \theta^q_\alpha(r) \sum_{k\neq\nul}  b_k^q \theta^q_\alpha(M^{*j}k) \\
	& \le
 	2^{\frac{2q}{p}}  C_{\w\phi}^q \|M^{*-j} \|^{sq} \sum_{r\in D(M^{*j})} | r|^{sq}  |\h f(r ) |^q  \theta^q_\alpha(r)    \|M^{*j}\|^{\alpha q} 	 \|b\|^q_{\ell_q^\alpha}\\
 &\le
	2^{\frac{2q}{p}}  C_{\w\phi}^q   \|b\|^q_{\ell_q^\alpha} \|M^{*j}\|^{\alpha q}  \|M^{*-j} \|^{sq}  \|f\|^{q,In}_{A_q^{s+\alpha},j}.
\end{split}
\end{equation*}
 For the second term $J_{111}$, using Lemma~\ref{lem1} with $\alpha=0$,  the Strang-Fix conditions~\eqref{fPhiSF} and~\eqref{fWPhi},  inequalities~(\ref{103}) and  (\ref{103++}), we get
\begin{equation*}
\begin{split}
J_{111} &:= 2^{\frac{2q}{p}} \sum_{r\in D(M^{*j})} \sum_{k\neq\nul}  |  \h{\phi_j}(r + M^{*j}k) |^q
    \left| \sum_{\underset{l\neq -k}{l\neq\nul}} \h f(r + M^{*j}(k+l))
	\overline{\h{\w\phi_j}(r + M^{*j}(k+l))} \right|^q \theta^q_\alpha(r + M^{*j}k) \\
	& \le
  2^{\frac{2q}{p}}\sum_{r\in D(M^{*j})} \sum_{k\neq\nul}  b_k^q  |M^{*-j} r|^{sq}  \|M^{*j}\|^{ \alpha q}\theta^q_\alpha(M^{*-j}r) \theta^q_\alpha(k)\\
    &\qquad\qquad\qquad\qquad\qquad\qquad\qquad\,\times\left| \sum_{\underset{l\neq -k}{l\neq\nul}} \h f(r + M^{*j}(k+l))
	\overline{\h{\w\phi_j}(r + M^{*j}(k+l))} \right|^q \\
	& \le
 2^{\frac{2q}{p}}   C_d^{q(s+\alpha)} \|M^{*j}\|^{ \alpha q}
\sum_{k\neq\nul}   b_k^q  \theta^q_\alpha(k)\sum_{r\in D(M^{*j})}  \left| \sum_{\underset{l\neq -k}{l\neq\nul}} \h f(r + M^{*j}(k+l))
	\overline{\h{\w\phi_j}(r + M^{*j}(k+l))} \right|^q \\
	&\le
 2^{\frac{2q}{p}}   C_d^{q(s+\alpha)}   \|b\|^q_{\ell_q^\alpha}	C^q_{Lem} \|M^{*j}\|^{ \alpha q} \|M^{*-j}\|^{q \gamma}  
	\|f\|^{q,Out}_{A_q^{\gamma},j}.
\end{split}
\end{equation*}
Thus, combining the estimates for $J_{11}$, $J_{110}$, and $J_{111}$, we derive
\begin{equation*}
  \begin{split}
     J_{11} \le &2^{\frac{2q}{p}}  C_{\w\phi}^q  \|b\|^q_{\ell_q^\alpha}  \|{M^{*j}}\|^{q\alpha}
\|M^{*-j} \|^{sq}   \|f\|^{q,In}_{A_q^{s+\alpha},j}\\
&+
2^{\frac{2q}{p}}
  C_d^{q(s+\alpha)} \|b\|^q_{\ell_q^\alpha} C^q_{Lem} \|{M^{*j}}\|^{ q\alpha} \|M^{*-j}\|^{q\gamma} \|f\|^{q,Out}_{A_q^\gamma,j}.
   \end{split}
\end{equation*}
Finally, the estimates for $J_1$, $J_{10}$, and $J_{11}$ yield
\begin{equation*}
\begin{split}
J_1 &\le C_{J_{10}} \|{M^{*j}}\|^{q\alpha} \|M^{*-j}\|^{q\gamma} \|f\|^{q,Out}_{A_q^\gamma,j} \\
 & \hspace{1cm}+
C_{J_{110}} \|{M^{*j}}\|^{q\alpha} \|M^{*-j} \|^{sq}  \|f\|^{q,In}_{A_q^{s+\alpha},j} + C_{J_{111}} \|{M^{*j}}\|^{ q\alpha} \|M^{*-j}\|^{q\gamma} \|f\|^{q,Out}_{A_q^\gamma,j},
\end{split}
\end{equation*}
which together with~\eqref{J0} implies
\begin{equation*}
\begin{split}
	&\sum_{n\in\zd} \left| \h f(n) - \sum_{l\in\zd} \h f(n + M^{*j}l)
	\overline{\h{\w\phi_j}(n + M^{*j}l)} \h {\phi_j}(n) \right|^q 	\theta^q_\alpha(n) \\
	&\le
2^{\frac qp} b_0^q  \|M^{*-j}\|^{q s} \|f\|^{q,In}_{A_q^{s+\alpha},j} + 2^{\frac qp} C^q_{\phi} C^q_{Lem}\|{M^{*j}}\|^{q\alpha} \|M^{*-j}\|^{q\gamma}\|f\|^{q,Out}_{A_q^\gamma,j} \\
&\hspace{1cm}+
C_{J_{110}}  \|{M^{*j}}\|^{q\alpha} \|M^{*-j} \|^{sq}  \|f\|^{q,In}_{A_q^{s+\alpha},j} + (C_{J_{10}} + C_{J_{111}} )\|{M^{*j}}\|^{ q\alpha} \|M^{*-j}\|^{q\gamma} \|f\|^{q,Out}_{A_q^\gamma,j} \\
&\le
C_1   \|{M^{*j}}\|^{q\alpha} \|M^{*-j}\|^{q s}  \|f\|^{q,In}_{A_q^{s+\alpha},j} + C_2 \|{M^{*j}}\|^{ q\alpha} \|M^{*-j}\|^{q\gamma}\|f\|^{q,Out}_{A_q^\gamma,j},
\end{split}
\end{equation*}
which completes the proof of~\eqref{lemMain} for the case $q<\infty$.

If now $q=\infty$, then
repeating step by step  all above estimates, one can easily obtain~\eqref{lemMain} for this case. \ $\Diamond$

\begin{rem}
Analyzing the proof of  Lemma~\ref{lem2}, it is not difficult see that
 $$C_1 = 2^{\frac{q}{p}}  b_0^q +  2^{\frac{2q}{p}}  C_{\w\phi}^q   \|b\|^q_{\ell_q^\alpha}, \quad
C_2\le C_2'\|b\|_{\ell_\infty} +C_2''\|b\|^q_{\ell_q^\alpha},
$$
where $C_2'$ and $C_2''$ do not depend on $b$.
\end{rem}

\begin{rem}
If under the assumptions of Lemma~\ref{lem2}  the functions $\phi_j$ are trigonometric polynomials such that ${\rm spec\ } \phi_j\subset {M^*}^jB_R(\nul)$, then
the assumption $\{ b_n\}_n \in \ell_q^\alpha$  is not required. Indeed,  there exists a finite set $\Omega\subset\zd$,
depending only  on $R$ and $d$,  such that $\h{\phi_j}({M^*}^jn+r)\ne0$
for at least one $r\in D(M^{*j})$ only for $n\in\Omega$. Analyzing the proof of Lemma~\ref{lem2}, we see that
 $\|b\|_{\ell_q^\alpha}$ can be replaced by $ C_{R,\alpha}\|b\|_{\ell_\infty}$.
\label{rem1}
\end{rem}

\begin{theo}
\label{theoMain}
Let   $1\le q \le \infty$, $1/p+1/q=1$, $j\in \N$, $\alpha\ge 0$, and let
	$\{\phi_j\}_j$, $\{\w\phi_j\}_j$ be a sequence of periodic functions in $L_p$ and periodic distributions, respectively,
	whose Fourier coefficients are such that conditions~(\ref{fWPhi}), (\ref{fPhiSF}), (\ref{fPhiPhiW}), (\ref{fPhi}) are valid with
	parameters $N\ge 0$, $s>0$, and $\{b_n\}_n \in \ell^\alpha_q$.
 Suppose $f\in A_q^{\gamma}$, where $\gamma>\alpha$ and $\gamma>N+ d/p$ for $q\ne1$ or  $\gamma\ge N$ for $q=1$.
Then
\be
\label{25}
\left\|
f - Q_j(f,\phi_j,\w\phi_j)
\right\|_{A_q^{\alpha}} \le
 C \|{M^{*j}}\|^{2\alpha} \|M^{*-j}\|^{\min\{\gamma,s+\alpha\}}\|f\|_{A_q^{\gamma}},
\ee
where   $C$ does not depend on $j$ and $f$.

In addition, if  $M$ is an isotropic matrix and $\lambda$
is its eigenvalue, then
\be
\label{27}
\left\|
f - Q_j(f,\phi_j,\w\phi_j)
\right\|_{A_q^{\alpha}}\le C'
|\lambda|^{-j\min\{\gamma-2\alpha,s-\alpha\}}\|f\|_{A_q^{\gamma}},
\ee
where   $C'$ does not depend on $j$ and $f$.
\end{theo}
{\bf Proof.}
First of all we mention that the inner product  $\langle f, \w\phi_j(\cdot - M^{-j} k)\rangle $ has meaning
under our assumptions, and hence the operator $Q_j(f,\phi_j,\w\phi_j)$ is well defined.
Taking into account that
\begin{equation*}
 \sum\limits_{k\in D(M^{j})}  e^{2 \pi i (k,M^{*-j}r)} =
\begin{cases}
m^j, &\mbox{if } \quad  r \equiv \nul\  (\text{mod} M^{*j}),
\\
0, &\mbox{if } \quad  r \not\equiv \nul\  (\text{mod} M^{*j}),
\end{cases}
\end{equation*}
we derive the following representation for the Fourier coefficients of $g:=f- Q_j(f,\phi_j,\w\phi_j)$:
\begin{equation*}
\begin{split}
\h g(n)
&=
\h f(n) - \frac{1}{m^j} \sum\limits_{k\in D(M^{j})} \langle f, \w\phi_j(\cdot - M^{-j} k)\rangle \h{\phi_j}(n) e^{-2\pi i (n, M^{-j}k)}
 \\
&=
\h f(n) - \frac{1}{m^j} \h{\phi_j}(n) \sum_{l\in\zd}
 \h f(l) \overline{\h{\w\phi_j}(l)}  \sum\limits_{k\in D(M^{j})}  e^{2 \pi i (k,M^{*-j}(l-n))}
 \\
&=
\h f(n) - \h{\phi_j}(n) \sum_{l\in\zd}
 \h f(n+M^{*j}l) \overline{\h{\w\phi_j}(n+M^{*j}l)}.
\end{split}
\end{equation*}
Using this  together with  Lemma~\ref{lem2}, taking into account
that ${\|M^{*j}\|^{-1}} \le \|M^{*-j}\|$, we obtain
\be
\label{25'}
\left\|
f - Q_j(f,\phi_j,\w\phi_j)
\right\|_{A_q^{\alpha}} \le
 C_1 \|{M^{*j}}\|^{2\alpha} \|M^{*-j}\|^{s+\alpha}\|f\|_{A_q^{s+\alpha}}+
C_2\|{M^{*j}}\|^{\alpha} \|M^{*-j}\|^{\gamma}
\|f\|_{A_q^{\gamma}}.
\ee
Assume first that $s+\alpha\le\gamma$.
 Obviously, in this case, $\|f\|_{A_q^{s+\alpha}}\le\|f\|_{A_q^{\gamma}}$ and
$\|M^{*-j}\|^{\gamma} \le  \|M^{*-j}\|^{s+\alpha}$ for all $j$ such that $\|M^{*-j}\| \le 1$. By~\eqref{00}, the last inequality is valid for all $j$ which are greater than some appropriate $j_0$, and there exists a constant $C_{j_0}$ such that $\|M^{*-j}\|^{\gamma} \le C_{j_0} \|M^{*-j}\|^{s+\alpha}$ for all $j \le j_0$, which yields~(\ref{25}).

Next, if $s+\alpha> \gamma$,  then we set $s'=\gamma-\alpha$ and note that
all assumptions of Theorem~\ref{theoMain} with $s'$ instead of $s$  are satisfied. Indeed, we need to check that~\eqref{fPhiSF} and~\eqref{fPhiPhiW} are valid for $s'$. For low dimensions ($d\le 4$), we have $|M^{*-j}r|\le 1$ for $r\in D(M^{*j})$ and, therefore, \eqref{fPhiSF} and~\eqref{fPhiPhiW} obviously hold. For $d>4$,~\eqref{fPhiSF} and~\eqref{fPhiPhiW} are valid for $s'$ if we replace the constants $b_n$ by $b_n \left(\frac{\sqrt{d}}{2}\right)^{s-s'}.$
Hence, inequality~\eqref{25'} with $s'$ instead of $s$ holds and $s'+\alpha=\gamma=\min\{s+\alpha, \gamma\}$, which yields~\eqref{25}.

Finally, inequality (\ref{27}) follows immediately from~(\ref{25}) and relation~(\ref{10}).  \ $\Diamond$

\vspace{.3cm}

Note that estimate~\eqref{27} is actually a generalization of~\eqref{spr}.
To compare these relations, one has to take into account that~\eqref{spr} was obtained under the following  strengthened  Strang-Fix conditions of order $s$ for a sequence of fundamental interpolants $\{\Lambda_j\}_j$:
$$
	| \h{\Lambda_j}(2^{j}n+r)| \le b_n |r|^s 2^{-j(s+\alpha)}, \quad \forall n\neq\nul,\, r\in D(M^{*j}).
$$
Analyzing the proof of Lemma~\ref{lem2}, it is not difficult to see that under the same strengthened  Strang-Fix conditions
for $\{\phi_j\}_j$ instead of~\eqref{fPhiSF}, we can replace $\min\{\gamma-2\alpha, s-\alpha\}$ by $\min\{\gamma-\alpha, s\} $ in inequality~\eqref{27}.

\begin{rem}
If under the assumptions of Theorem~\ref{theoMain}, the functions $\phi_j$ are trigonometric polynomials such that ${\rm spec\ } \phi_j\subset D(M^{*j})$, then, analyzing the proof of   Lemma~\ref{lem2}, one can easily see that
$$J_1 = J_{100} = \sum_{n\notin D(M^{*j})}  | \h f(n) |^q  \theta^q_\alpha(n) \le
 2^{q (\gamma-\alpha)} \|M^{*-j}\|^{q(\gamma-\alpha)} \|f\|^{q,Out}_{A_q^{\gamma},j},
$$
which together with~\eqref{J0} implies
  $$
  J \le
  C_1  \|f\|^{q,In}_{A_q^{s+\alpha},j}  \|M^{*-j}\|^{q s} + C_2 \  \
	\|{M^{*j}}\|^{q\alpha} \|M^{*-j}\|^{q\gamma} \|f\|^{q,Out}_{A_q^\gamma,j}.
  $$
	It follows that (\ref{25}) can be replaced by
\begin{equation*}
\left\|
f - Q_j(f,\phi_j,\w\phi_j)
\right\|_{A_q^{\alpha}} \le
C \|{M^{*j}}\|^{\alpha} \|M^{*-j}\|^{\min\{\gamma,s+\alpha\}}\|f\|_{A_q^{\gamma}},
\end{equation*}
where   $C$ does not depend on $j$ and $f$.
In a similar way, relation~(\ref{27})  can be also improved in this case.
\label{rem2}
\end{rem}

Next we need the following embedding properties  between the spaces $L_p$ and $A_q^0$.
\begin{prop}
If $2 \le p \le \infty$, $1/p+1/q=1$,  then $A_q^{0}  \subset L_p$,  with
\begin{equation}
\|f\|_p \le \|f\|_{A_q^0}.
\label{fT12}
\end{equation}
\label{T1}
\end{prop}

Relation~\eqref{fT12} is the classical Hausdorff-Young inequality, for its multivariate version  see, e.g.,~\cite[p.~174]{Gr} or~\cite{Spr2000}.

\begin{theo}
\label{big_p_th}
Let  $2\le p\le\infty$,  $1/p+1/q=1$,  $\{\phi_j\}_j$, $\{\w\phi_j\}_j$,  $N$,  $s$, and $\gamma$ be as in Theorem~\ref{theoMain} with
$\alpha=0$ and   $f\in A_q^\gamma$. Then
\be
\label{125}
\left\|
f - Q_j(f,\phi_j,\w\phi_j)
\right\|_p \le
C_1     \|M^{*-j}\|^{s}\|f\|^{In}_{A_q^{s},j} + C_2 \   \|M^{*-j}\|^{\gamma} \|f\|^{Out}_{A_q^\gamma,j},
\ee
where  $C_1$ and $C_2$ do not depend on $j$ and $f$.

If $\vartheta$ is any positive number  which is smaller in absolute value than any eigenvalue of $M$, then
\be
\label{130}
\left\|
f - Q_j(f,\phi_j,\w\phi_j)
\right\|_p
\le C\vartheta^{-j\min\{s,\gamma\}}\|f\|_{A_q^{\gamma}},
\ee
where $C$ does not depend on $j$ and $f$.
In addition, if  $M$ is an isotropic matrix and $\lambda$
is its eigenvalue, then $\vartheta$  can be replaced by $|\lambda|$ in~\eqref{130}.
\end{theo}

{\bf Proof.}
Obviously, the operator $Q_j(f,\phi_j,\w\phi_j)$ is well defined and belongs to $L_p$ under our assumptions.
Using~\eqref{fT12}, one can see that $f\in L_p$ and
\be
\label{131}
\left\|
f - Q_j(f,\phi_j,\w\phi_j)
\right\|_p \le
\left\|
f - Q_j(f,\phi_j,\w\phi_j)
\right\|_{A_q^0}.
\ee
Analyzing the proof of~Theorem~\ref{theoMain} and using Lemma~\ref{lem2},
we obtain~\eqref{125}.
Inequality~\eqref{130} follows from~\eqref{131}, ~\eqref{25} with $\alpha=0$, and~\eqref{00}.
It remains to say that in the case of isoropic matrix $M$, one can use~\eqref{10} instead of~\eqref{00}.
   \ $\Diamond$

Note that non-periodic  counterparts of the results stated in  Theorem~\ref{big_p_th} were obtained in~\cite{KS}, see Theorems~4, 5.

\begin{coro}
	\label{big_p}
	If under the assumptions of Theorem~\ref{big_p_th} we have
	$\h f(n)= \mathcal{O}(|n|^{-\kappa})$, where $\kappa>N+d$, and $M$ is an isotropic  matrix, then
	\be
	\label{127}
	\left\|
	f - Q_j(f,\phi_j,\w\phi_j)
	\right\|_p =
	\begin{cases}
		\mathcal{O}\big(|\lambda|^{-js}\big), & s< \kappa-d/q,
		\\
		\mathcal{O}\big(j |\lambda|^{-js}\big), & s= \kappa-d/q,
		\\
		\mathcal{O}\big(|\lambda|^{-j(\kappa-\frac dq)}\big), & s>\kappa-d/q,
		
	\end{cases}
	\ee	
	where  $\lambda$ is an eigenvalue of $M$.
\end{coro}

{\bf Proof.}
First, we assume that $s< \kappa - \frac dq$ and choose $\gamma>N+\frac dp$ such that $s<\gamma <  \kappa - \frac dq$. Therefore, $s-\kappa < \gamma-\kappa < - \frac dq$ and both expressions $\|f\|^{q,Out}_{A_q^\gamma,j}$ and $\|f\|^{q,In}_{A_q^s,j},$ in~\eqref{125} are finite. Hence, inequality~\eqref{125} together with~\eqref{10} yields that
$$
\left\|
f - Q_j(f,\phi_j,\w\phi_j)
\right\|^q_p \le C \|M^{*-j}\|^{qs} \sum_{k\in\zd} |k|^{(\gamma-\kappa)q} = \mathcal{O}(\|M^{*-j}\|^{qs}).
$$

Now, we assume that $s\ge \kappa-\frac dq$. Choose $\gamma>N+\frac dp$ such that $\gamma <  \kappa - \frac dq$. Let us consider the second term in~\eqref{125}. Since  $\frac{1}{|n|}\le \frac{\|M^{*-j}\|}{|M^{*-j} n|}$ we get
\begin{equation*}
\begin{split}
\|M^{*-j}\|^{q\gamma} \|f\|^{q,Out}_{A_q^\gamma,j} & \le
\|M^{*-j}\|^{q\gamma} \sum_{n\neq \nul} \sum_{r \in D(M^{*j})} \frac{1}{|M^{*j}n + r|^{(\kappa-\gamma)q}} \\
&\le
\|M^{*-j}\|^{q\gamma}   \sum_{n\neq \nul}   \sum_{r \in D(M^{*j})} \max_{\xi \in [-\frac 12,\frac 12)^d} \frac{\|M^{*-j}\|^{q(\kappa - \gamma)}}{|n +\xi|^{(\kappa-\gamma)q}}\\
&\le C m^j \|M^{*-j}\|^{q\kappa}
\end{split}
\end{equation*}
Next, consider the first term in~\eqref{125}. Then
$$
 \|M^{*-j}\|^{q s}\|f\|^{q,In}_{A_q^{s},j} \le
C \|M^{*-j}\|^{q s}  \sum_{n \in D(M^{*j})}  |n|^{(s-\kappa)q}.
$$
It is clear that there exists $A = A(d) > \sqrt{d}/2$ such that for $|n| \ge A$, $|x|/2\le |n| \le 2|x|$ for any $x\in n + \frac 12[-1,1)^d$ (note that anyway $|x|> \w A$, for some $\w A>0$). Therefore, if $s - \kappa \ge 0$, then $ |n|^{(s-\kappa)q} \le 2^{(s-\kappa)q} |x|^{(s-\kappa)q}$, and if  $s - \kappa <  0$,  then $|n|^{(s-\kappa)q} \le  (|x|/2)^{(s-\kappa)q}$. Thus,
 $$|n|^{(s-\kappa)q} \le 2^{|s-\kappa|q} \int_{n+ \frac 12[-1,1)^d} |x|^{(s-\kappa)q} dx.$$
Hence,
\begin{equation*}
\begin{split}
\sum_{n \in D(M^{*j})}  |n|^{(s-\kappa)q} & = \sum_{|n| < A}   |n|^{(s-\kappa)q}  + \sum_{|n| \ge A}   |n|^{(s-\kappa)q} \le
\w C_1 +  2^{|s-\kappa|q}  \int\limits_{\w A<|x|<\|M^{*j}\|\sqrt{\frac d4}}  |x|^{(s-\kappa)q} dx \\
& \le
\w C_2 \int\limits_{\w A}^{\|M^{*j}\|\sqrt{\frac d4}} r^{(s-\kappa)q  +d - 1} dr  \le
\begin{cases}
\w C_3 \log \|M^{*j} \|,\  s= \kappa - \frac dq;\\
\w C_4 \|M^{*j}\|^{(s-(\kappa-\frac dq))q },\  s > \kappa-\frac dq.
\end{cases}
\end{split}
\end{equation*}
Overall, for $ s= \kappa - \frac dq$ inequality~\eqref{125} and above considerations yields that
$$
\left\|
f - Q_j(f,\phi_j,\w\phi_j)
\right\|^q_p  \le  C \|M^{*-j}\|^{q s} (m^j \|M^{*-j}\|^{d}  +   \log \|M^{*j} \|).
$$
For  $ s>\kappa - \frac dq$ inequality~\eqref{125} and above considerations yields that
$$
\left\|
f - Q_j(f,\phi_j,\w\phi_j)
\right\|^q_p  \le  C ( m^j \|M^{*-j}\|^{q \kappa}   +    \|M^{*-j}\|^{q s} \|M^{*j}\|^{(s-(\kappa-\frac dq))q}).
$$
It remains to note that $m=|\lambda|^d$ for isotropic dilation matrix, and using~\eqref{10},  we get the required estimates. $\Diamond$

\vspace{.5cm}

Let   $\mathcal{B}=\mathcal{B}_{\delta, R}$ denote the class of sequences of trigonometric polynomials $\phi_j$, $j\in\n$, such that
\begin{equation}
\label{my1}
{\rm spec\ } \phi_j\subset {M^*}^j\Big(B_R(\nul)\setminus\bigcup_{n\in\zd\setminus\{\nul\}} B_\delta(n)\Big),
\end{equation}
\begin{equation}
\label{my2}
|\h{\phi_j}(l)|\le C'_\phi\quad \forall j\in\n, \forall l\in\zd,
\end{equation}
for some positive  constants  $R, \delta$, and $C'_\phi$.

Note that  $\{\phi_j\}_{j}$ belongs to $ \mathcal{B} $ whenever (\ref{my2}) is valid and
$$
{\rm spec\ } \phi_j\subset M^{*j}\big([-1+\delta, 1-\delta]^d\big),\quad \delta \in (0,1/2).
$$

\begin{defi}
We say that a sequence  of integrable periodic functions $\{\phi_j\}_{j}$ is strictly compatible with a sequence of periodic distributions
 $\{\w\phi_j\}_{j}$ if
\begin{equation}
\label{my3}
\h{\w\phi_j}(l)\overline{\h{\phi_j}(l)}=1,\quad \forall j\in\n, \forall l\in\zd:\ |M^{*-j}l|\le \delta,
\end{equation}
for some $\delta>0$.
\end{defi}

\begin{theo}
\label{theo_my1}
Let   $1\le q < \infty$, $1/p+1/q=1$,
  $\{\w\phi_j\}_{j}$ be a sequence of periodic distributions satisfying~(\ref{fWPhi}) with $N\ge0$,
and a sequence of trigonometric polynomials $\{\phi_j\}_{j}\in \mathcal{B}$ be  strictly compatible with  $\{\w\phi_j\}_{j}$
with respect to the parameter $\delta\in (0,1/2)$.
Suppose $f\in A_q^{\gamma}$, where $\gamma>0$  and $\gamma>N+ d/p$ for $q\ne1$ or  $\gamma\ge N$ for $q=1$,
and $\alpha\in [0,\gamma)$.
Then
\be
\label{my4}
\left\|
f - Q_j(f,\phi_j,\w\phi_j)
\right\|_{A_q^{\alpha}} \le C
\|{M^{*j}}\|^{2\alpha} \|{M^*}^{-j}\|^{\gamma}
\(\sum\limits_{|M^{*-j}r|\ge \delta}|r|^{q\gamma}|\h f(r)|^{q}\)^\frac1q,
\ee
and moreover, in the case  $p\ge2$,
\be
\label{my5}
\left\|
f - Q_j(f,\phi_j,\w\phi_j)
\right\|_p\le C'\|{M^*}^{-j}\|^{\gamma}
\(\sum\limits_{|M^{*-j}r|\ge \delta}|r|^{q\gamma}|\h f(r)|^{q}\)^\frac1q,
\ee
where   $C$, $C'$ do not depend on $f$ and $j$.
\end{theo}
{\bf Proof.}
As in Theorem~\ref{theoMain}, the operator $Q_j(f,\phi_j,\w\phi_j)$ is well defined under our assumptions.
Let us check that for any $s>0$ the Strang-Fix conditions of order $s$ for $\phi_j$
and the weak compatibility conditions for $\phi_j$ and $\w \phi_j$ of order $s$ are satisfied.
Indeed, let $n\in\zd, n\neq\nul,\, r\in D(M^{*j})$. If $|M^{*-j} r|\le \delta$, then $\h{\phi_j}(M^{*j}n+r)=0$;
if $|M^{*-j} r|> \delta$, then, by~(\ref{my2}),
$$
| \h{\phi_j}(M^{*j}n+r)| \le C'_\phi \delta^{-s}\delta^{s}\le C'_\phi \delta^{-s} |M^{*-j} r|^s.
$$
Thus, conditions (\ref{fPhiSF}) are satisfied with  $b_n\equiv C'_\phi \delta^{-s} $.
 Similarly, for every $ r\in D(M^{*j})$, using~(\ref{my2}), (\ref{my3})  and~(\ref{fWPhi}), we have
$$
|1 - \overline{\h{\w\phi_j}(r)}  \h{\phi_j}(r)| \le(1+C'_\phi C_{\w\phi})
\delta^{-s}\delta^{s}\le (1+C'_\phi C_{\w\phi})\delta^{-s} |M^{*-j} r|^s,
$$
which implies that (\ref{fPhiPhiW}) is satisfied with $b_0=(1+C'_\phi C_{\w\phi})\delta^{-s}$.

Now, taking into account Remark~\ref{rem1}, we see that all assumptions of
 Lemma~\ref{lem2}  are fulfilled with any $s>0$.
 Analyzing the proof of this lemma, we see that
\begin{equation}
\begin{split}
	J & := \sum_{n\in\zd} \left| \h f(n) - \h {\phi_j(n)} \sum_{l\in\zd} \h f(n + M^{*j}l)
	\overline{\h{\w\phi_j}(n + M^{*j}l)}\right|^q \theta^q_\alpha(n)
	\\
	&\le
 J_{00}+J_{110} + C_2  \|M^{*j}\|^{\alpha q} \|M^{*-j}\|^{q\gamma} \|f\|^{q,Out}_{A_q^\gamma,j},
\label{101}
\end{split}
\end{equation}
Obviously,
$$
\|f\|^{q,Out}_{A_q^\gamma,j}\le   2^{q\gamma}
\sum\limits_{|M^{*-j}r|\ge  \delta}|r|^{q\gamma}|\h f(r)|^{q\gamma}.
$$

Repeating step by step the estimates for $J_{00}$ and $J_{110}$, using~(\ref{my3}), (\ref{my1}), and Remark~\ref{rem1},
and taking into account that ${\|M^{*j}\|^{-1}} \le \|M^{*-j}\|$, we obtain
\begin{equation}
\begin{split}
J_{00} &= 2^{q/p} \sum_{r\in D(M^{*j})} |\h f(r)|^q | 1 -
	  \overline{\h{\w\phi_j}(r)}\h{\phi_j}(r)|^q \theta^q_\alpha(r)  \\
	  &\le
2^{q(\alpha+1/p)} b_0^q \|M^{*j}\|^{\alpha q} \|M^{*-j} \|^{(s+\alpha)q}
\sum\limits_{|M^{*-j}r|\ge \delta}|r|^{q(s+\alpha)}|\h f(r)|^{q}
\end{split}
\end{equation}
and
\begin{equation*}
\begin{split}
J_{110} & =  2^{\frac{2q}{p}}  \sum_{r\in D(M^{*j})} \sum_{k\neq\nul}  | \h{\phi_j}(r + M^{*j}k) |^q
  |\h f(r )  	\overline{\h{\w\phi_j}(r )}|^q \theta^q_\alpha(r + M^{*j}k) \\
  &\le
	  2^{q(\alpha+2/p)}\delta^{-s} C_\phi' C_{\w\phi}^q  C_{R,\alpha}^q \|M^{*j}\|^{2\alpha q} \|M^{*-j} \|^{(s+\alpha)q}
	\sum\limits_{|M^{*-j}r|\ge \delta}|r|^{q(s+\alpha)}|\h f(r)|^{q},
\end{split}
\end{equation*}
where $C_{R,\alpha}$ is the constant from Remark~\ref{rem1}.
Substituting these relations into~(\ref{101}) and choosing $s=\gamma-\alpha$, we obtain~\eqref{my4}.

To prove~\eqref{my5} it remains to set $\alpha=0$ and combine~\eqref{my4}  with~\eqref{fT12}.  \   $\Diamond$
\vspace{.3cm}

Note that  relation~\eqref{my5} with $p=\infty$ is a periodic analog of Theorem~15 in~\cite{KKS}, which, in turn, is a generalization of Brown's inequality, see~\cite{Brown}.

\begin{prop}
\label{coro1}
If under the assumptions of Theorem~\ref{theo_my1}, equality~\eqref{my3} holds for all $l\in D(M^{*j})$
and $\spec\phi_j\subset D(M^{*j})$, then
\be
\label{my4'}
\left\|
f - Q_j(f,\phi_j,\w\phi_j)
\right\|_{A_q^{\alpha}} \le C
\|{M^{*j}}\|^{\alpha}  \|{M^*}^{-j}\|^{\gamma}
\(\sum\limits_{r\not\in D(M^{*j})}|r|^{q\gamma}|\h f(r)|^{q}\)^\frac1q,
\ee
and moreover, in the case  $p\ge2$,
\be
\label{my5'}
\left\|
f - Q_j(f,\phi_j,\w\phi_j)
\right\|_p\le C'\|{M^*}^{-j}\|^{\gamma}
\(\sum\limits_{r\not\in D(M^{*j})}|r|^{q\gamma}|\h f(r)|^{q}\)^\frac1q,
\ee
where   $C$ and $C'$ do not depend on $f$ and $j$.
\end{prop}

{\bf Proof.} Analyzing the proof  of Theorem~\ref{theo_my1},  we see that under our assumptions $J_{00}=J_{110}=0$ and
$$
\|f\|^{q,Out}_{A_q^\gamma,j}\le   2^{q\gamma}
\sum\limits_{r\not\in D(M^{*j})}|r|^{q\gamma}|\h f(r)|^{q\gamma}.
$$
Using this and taking into account Remark~\ref{rem2}, we obtain~\eqref{my4'} and \eqref{my5'}.\ $\Diamond$

\subsection{$L_p$-errors  using the best approximation}

In this section, we will show that in some partial cases the results obtained in the previous sections can be sharpened and extended to a wider class of functions $f$ than those considered in Theorems~\ref{theoMain}, \ref{big_p_th}, and~\ref{theo_my1}. Moreover, we show that the $L_p$-error estimates for $Q_j(f,\phi_j,\w\phi_j)$ can be given using the best approximation.

We restrict ourselves to the case of a diagonal dilation matrix
$M={\rm diag} (m_1, m_2,\dots,m_d)$, $m_j\in \Z$,
and the case $\{\vp_j\}_j\in \mathcal{B}$, where $\mathcal{B}=\mathcal{B}_{\delta, R}$ is the class of sequences of
trigonometric polynomials introduced in Section~4.1.

We need to specify the class of tempered distributions $\{\w\phi_j\}_j$. We will say that a sequence of tempered distributions $\{\w\phi_j\}_j$ belongs to the class $\mathcal{S}_{N,p}'$ for some $N\ge 0$ and $1\le p\le \infty$ if it satisfies~\eqref{fWPhi} and for any trigonometric polynomial $T_n$ such that $\spec T_n\subset \{k\in \Z^d\,:\, |M^{-n}k|\le 1\}$, one has
\begin{equation}\label{DefS}
  \Vert T_n * \w\phi_j\Vert_p \le C_{\w\vp,p} m^{\frac Nd (n-j)} \Vert T_n \Vert_p\quad\text{for all}\quad n\ge j,\quad j,n\in \N.
\end{equation}

Note that in the case of $p=2$ and an isotropic matrix $M$, conditions~\eqref{fWPhi} imply inequality~\eqref{DefS}.
As a simple example of $\{\w\phi_j\}_j \in \mathcal{S}_{N,p}'$, we can take distributions $\{\w\phi_j\}_j$ corresponding to some differential operator. Namely, let $d=1$ and
$\h{\w\phi_j} (l) = \sum_{0\le\beta\le N} c_\beta (2\pi i M^{-j} l)^{\beta}$, $N\in \z_+$, $j\in \N$,
then by the well-known Bernstein inequality for trigonometric polynomials
$\Vert T_n^{(r)}\Vert_{L_p(\T)} \le n^r \Vert T_n \Vert_{L_p(\T)}$
(see, e.g.,~\cite[p.~215]{Timan}),
we easily derive that $\{\w\phi_j\}_j \in \mathcal{S}_{N,p}'$.

We will also use the following class of sequences $\mathcal{L}_{p}$, $1\le p\le \infty$. We will say that a sequence
$\{\w\phi_j\}_j$ belongs to $\mathcal{L}_{p}$ if for all $j\in \N$
$$
\Vert \w\phi_j\Vert_{\mathcal{L}_{p,j}}:=\(m^j\int_{M^{-j}\T^d} \bigg(\frac1{m^j}\sum_{k\in D(M^j)} |\w\phi_j(x-M^{-j}k)| \bigg)^p dx\)^\frac1p<\infty\quad\text{if}\quad 1\le p<\infty
$$
and
$$
\Vert \w\phi_j\Vert_{\mathcal{L}_{\infty,j}}:=\frac1{m^j}\sup_{x\in \R^d}\sum_{k\in D(M^j)} |\w\phi_j(x-M^{-j}k)| <\infty\quad\text{if}\quad p=\infty.
$$

An important example of $\{\w\phi_j\}_j\in \mathcal{L}_{p}$ is given by the normalized characteristic functions of $M^{-j} [-\frac12,\frac12)^d$,
$
\w\phi_j (x)=m^j \chi_{M^{-j} [-\frac12,\frac12)^d}(x).
$
It is worth noting that such $\w\phi_j$ provide a periodic counterpart of Kantorovich-Kotelnikov operators studied in~\cite{KS3}.

For any $d\times d$-matrix $A$, we denote
$$
\mathcal{T}_{A}:=\{T\,:\, \spec T\subset \{k\in \Z^d\,:\, |A^{-1}k|<1\}\}.
$$
The error of the best approximation of $f\in L_p$ by trigonometric polynomials $T\in \mathcal{T}_{A}$ is defined by
$$
E_{A}(f)_p:=\inf\left\{\Vert f-T\Vert_p\,:\, T\in \mathcal{T}_{A}\right\}.
$$

We will use the following anisotropic Besov spaces with respect to the matrix~$M$. We will  say that
$f\in \mathbb{B}_{p,q}^s (M)$, $1\le p\le\infty$, $0<q\le \infty$, and $s>0$, if $f\in L_p$ and
$$
\Vert f\Vert_{\mathbb{B}_{p,q}^s (M)}:=\Vert f\Vert_p+\(\sum_{\nu=1}^\infty m^{\frac sd q\nu} E_{M^\nu} (f)_p^q\)^{\frac 1q}<\infty.
$$

For simplicity, we will also denote
$$
\left\Vert \{a_k\}_{k}\right\Vert_{\ell_{p,M^j}}:=\left\{
                                                 \begin{array}{ll}
                                                  \displaystyle\bigg(\frac1{m^j}\sum\limits_{k\in D(M^j)}|a_k|^p\bigg)^\frac1p, & \hbox{if $1\le p<\infty$,} \\
                                                   \displaystyle\sup\limits_{k\in D(M^j)}|a_k|, & \hbox{if $p=\infty$.}
                                                 \end{array}
                                               \right.
$$

\begin{lem}\label{lemMZ}
Let $1\le p\le\infty$ and $T_j\in \mathcal{T}_{\lambda M^j}$ for some $\lambda>0$. Then, for any $j\in \N$,
$$
\left\Vert \{ T_j(M^{-j}k) \}_{k} \right\Vert_{\ell_{p,M^j}}\le C\Vert T_j\Vert_p,
$$
where the constant $C$ depends only on $p$, $d$, and $\lambda$.  
\end{lem}

{\bf Proof.} In the case $p=\infty$, the proof is obvious. The case $1\le p<\infty$ directly follows from
the following Marcinkiewicz-Zygmund type inequality (see~\cite{LMN})
$$
\sum_{j=1}^m |T_n(\tau_j)|^p\le (p+1)\frac e2 \(2n+\frac1\delta\)\int_0^1 |T_n(x)|^p dx,
$$
where $T_n$ is a univariate trigonometric polynomial of degree at most $n$, $-1/2\le \tau_1<\tau_2<\dots<\tau_m<1/2$, and
$\delta=\min\{\tau_2-\tau_1,\tau_3-\tau_2,\dots,\tau_m-\tau_{m-1},1-(\tau_m-\tau_1)\}.$~~$\Diamond$

\begin{lem}\label{lem1Be}
Let  $1\le p\le \infty$, $N\ge 0$, $\delta\in (0,1/2)$, $\w\vp\in \mathcal{S}_{N,p}'$,
and $f\in \mathbb{B}_{p,1}^{d/p+N}(M)$. Suppose that the polynomials
$T_\mu$, $\mu\in\z_+$,  are such that $\spec T_\mu \subset D(M^\mu)$ and
$$
\Vert f-T_\mu \Vert_p\le \sigma E_{\delta M^\mu}(f)_p
$$
with some constant $\sigma$ independent of $f$ and $\mu$. Then the sequence $\{\langle T_\mu, \w\phi_{j}(\cdot-M^{-j}k)\rangle\}_{\mu=0}^\infty$ converges uniformly with respect to $k\in\zd$ and $j\in \N$,  and the limit does not depend on the choice of polynomials~$T_\mu$. Moreover, for all $n\in \N$,
$$
	\sum_{\mu=n}^\infty \Vert \{\langle T_{\mu+1}-T_\mu, \w\phi_{j}(\cdot-M^{-j}k)\rangle \}_{k} \Vert_{\ell_{p,M^j}}\le
	 C m^{-(\frac 1p+\frac Nd)j}\sum_{\mu=n}^\infty m^{(\frac1p+\frac Nd)\mu} E_{\delta M^\mu}(f)_p,
	$$
where the constant $C$ depends only on $d$, $p$, $M$, and $\sigma$.
\end{lem}

{\bf Proof.} 
Let $n\ge j$, $n\in \N$. Using Lemma~\ref{lemMZ}, condition~\ref{DefS}, we obtain
\begin{equation}\label{1yk}
  \begin{split}
     &\sum_{\mu=n}^\infty \Vert \{\langle T_{\mu+1}, \w\phi_{j}(\cdot-M^{-j}k)\rangle -\langle T_\mu, \w\phi_{j}(\cdot-M^{-j}k)\rangle \}_{k} \Vert_{\ell_{p,M^j}}\\
     &=\sum_{\mu=n}^\infty \Vert \{\langle T_{\mu+1}-T_\mu, \w\phi_{j}(\cdot-M^{-j}k)\rangle\}_{k}\Vert_{\ell_{p,M^j}}=\sum_{\mu=n}^\infty \Vert \{ (T_{\mu+1}^--T_\mu^-)*\w\phi_{j} (M^{-j}k)\}_{k}\Vert_{\ell_{p,M^j}}\\
      &\le m^{-\frac jp}\sum_{\mu=n}^\infty m^{\frac \mu p} \Vert \{ (T_{\mu+1}^--T_\mu^-)*\w\phi_{j} (M^{-\mu}k)\}_{k}\Vert_{\ell_{p,M^\mu}}\le C_1 m^{-\frac jp}\sum_{\mu=n}^\infty m^{\frac \mu p} \Vert  (T_{\mu+1}^--T_\mu^-)*\w\phi_{j} \Vert_{p}\\
      &\le C_2 m^{-(\frac 1p+\frac Nd)j}\sum_{\mu=n}^\infty m^{(\frac1p+\frac Nd)\mu} \Vert T_{\mu+1}^--T_\mu^- \Vert_{p}
      \le C_2 m^{-(\frac 1p+\frac Nd)j}\sum_{\mu=n}^\infty m^{(\frac1p+\frac Nd)\mu} E_{\delta M^\mu}(f)_p.
  \end{split}
\end{equation}
The latter series is convergent since $f\in \mathbb{B}_{p,1}^{d/p+N}(M)$, which yields that
the sequence $\{\langle T_\mu, \w\phi_{j}(\cdot-M^{-j}k)\rangle\}_{\mu=1}^\infty$ is convergent for every $k\in\zd$ and $j\in \N$.
By analogy with~\eqref{1yk}, it is also easy to  check that the limit does not depend on
the choice of functions $T_\mu$.~~$\Diamond$

\medskip

In the previous section, the operator
$$
Q_j(f,\phi_j,\w\phi_j) = \frac{1}{m^j} \sum\limits_{k\in D(M^j)} \langle f, \w\phi_j(\cdot - M^{-j} k)
\rangle \phi_j (\cdot - M^{-j}k),
$$
was considered for functions $f$ whose Fourier coefficients decay sufficiently fast. In particular, supposing that \eqref{fWPhi} holds and  $f\in A_q^\gamma$ with $\gamma\ge N$, we defined the inner product by
\be
\label{100my+}
 \langle f, \w\phi_j(\cdot - M^{-j} k) \rangle
 =  \sum_{l\in\zd}
 \h f(l) \overline{\h{\w\phi_j}(l)}   e^{2 \pi i (k,M^{*-j}l)}.
\ee
Taking into account Lemma~\ref{lem1Be} and using condition~\eqref{DefS}, we extend the functional $\langle f, \w\phi_j(\cdot - M^{-j} k) \rangle$ from $f\in A_q^\gamma$ to the Besov spaces $\mathbb{B}_{p,1}^{d/p+N}(M)$ as follows.

	\begin{defi}
\label{def0}
	Let $1\le p\le \infty$, $N\ge 0$, $\w\phi\in S'_{N, p}$, and the polynomials $T_\mu$ be as in Lemma~\ref{lem1Be}.
	For every $f\in \mathbb{B}_{p,1}^{d/p+N}(M)$ and $k\in\zd$, we set
$$
\langle f,\w\phi_j(\cdot - M^{-j} k)\rangle:=\lim\limits_{\mu\to\infty}\langle T_\mu,\w\phi_j(\cdot - M^{-j} k)\rangle, \quad j\in\N.
$$
	\end{defi}	
	
Now, if $\w\phi\in S'_{N, p}$, then the quasi-projection  operators
$$
Q_j(f,\phi,\w\phi)=\sum_{k\in\zd}\langle f,\w\phi_j(\cdot - M^{-j} k)\rangle \phi_j(\cdot - M^{-j} k)
$$
are  defined on the space $\mathbb{B}_{p,1}^{d/p+N}(M)$ for a wide class of appropriate functions $\phi$.

\medskip

Note that below we will consider the operators $Q_j(f,\phi,\w\phi)$ for the sequences $\{\w\phi_j\}_j$ belonging to $\mathcal{L}_{q}$. In this case,
the inner product $\langle f,\w\phi_j(\cdot - M^{-j} k)\rangle$ has sense for any $f\in L_p$, $1/p+1/q=1$. As usual, we have $\langle f,\w\phi_j(\cdot - M^{-j} k)\rangle=\int_{\T^d} f(t)\overline{\w\phi_j}(t - M^{-j} k)dt$.

\medskip

In what follows, a Fourier multiplier operator associated with a function $\phi_j$ is denoted by $S_{\phi_j}$, i.e. for any function $f\in L_p$, $1\le p\le \infty$, we set 
$$
S_{\phi_j}(f;x):=\sum_{k\in \Z^d} \h{\phi_j}(k) \h f(k) e^{2\pi i(k,x)}.
$$
Denote also
$$
K_{\phi_j,q}:=\sup_{\Vert f\Vert_q\le 1}\Vert S_{\phi_j}(f;\cdot)\Vert_q.
$$

The standard example of such operators is the partial sums of Fourier series. For example, if
$\h{\phi_j}(\xi)=\chi_{M^j [-\frac12,\frac12)^d}(\xi)$, then $S_{\phi_j}$ represents the rectangular partial sums of Fourier series and
$$
K_{\phi_j,q} \asymp \left\{
                                 \begin{array}{ll}
                                   1, & \hbox{$1<q<\infty$,} \\
                                   j^d, & \hbox{$q=1$ or $\infty$.}
                                 \end{array}
                               \right.
$$

We will need the following two auxiliary lemmas.

\begin{lem}\label{lemKK1}
Let $1\le p\le \infty$, $1/p+1/q=1$, $\{a_k\}_{k\in D(M^j)}\in \mathbb{C}$, and $\{\phi_j\}_j \in \mathcal{B}$. Then, for any $j\in \N$,
\begin{equation*}
  \bigg\Vert \frac1{m^j}\sum_{k\in D(M^j)} a_k \phi_j(\cdot-M^{-j}k)\bigg\Vert_p
\le  C K_{\phi_j,q} \left\Vert \{a_k\}_{k}\right\Vert_{\ell_{p,M^j}}.
\end{equation*}
where the constant $C$ does not depend on $j$ and $\{a_k\}$.
\end{lem}
{\bf Proof.}
Consider the case $1\le p<\infty$. By duality, we can find a function  $g_j$ such that $\Vert g_j\Vert_q\le 1$ and
\begin{equation}\label{eqKK1}
\begin{split}
   \bigg\Vert \sum_{k\in D(M^j)} a_k \phi_j(\cdot-M^{-j}k)\bigg\Vert_p&=\bigg|\bigg\langle \sum_{k\in D(M^j)} a_k \phi_j(\cdot-M^{-j}k),g_j\bigg\rangle\bigg|\\
   &=\bigg|\sum_{k\in D(M^j)} a_k \langle \phi_j(\cdot-M^{-j}k),g_j \rangle\bigg|.
\end{split}
\end{equation}
Applying H\"older's inequality and Lemma~\ref{lemMZ}, taking into account that
$\langle \phi_j(\cdot-M^{-j}k),g_j \rangle=\phi_j*{g_j}^{-}(-M^{-j}k)$ and $\phi_j*g_j^-\in  \mathcal{T}_{R M^j}$,
we obtain
\begin{equation}\label{eqKK2}
\begin{split}
\sum_{k\in D(M^j)} &|a_k \langle \phi_j(\cdot-M^{-j}k),g_j \rangle|\le \bigg(\sum_{k\in D(M^j)} |a_k|^p \bigg)^\frac1p
\bigg(\sum_{k\in D(M^j)} |\langle \phi_j(\cdot-M^{-j}k),g_j \rangle|^q \bigg)^\frac1q\\
&\le C\bigg(\sum_{k\in D(M^j)} |a_k|^p \bigg)^\frac1p
m^{j/q}\Vert \phi_j*{g_j^-}\Vert_q \le C\bigg(\sum_{k\in D(M^j)} |a_k|^p \bigg)^\frac1p
m^{j/q}K_{\phi_j,q}\Vert g_j \Vert_q\\
&\le C m^j K_{\phi_j,q} \left\Vert \{a_k\}_{k}\right\Vert_{\ell_{p,M^j}}.
\end{split}
\end{equation}
Then, combining~\eqref{eqKK1} and~\eqref{eqKK2}, we prove the lemma for $p\neq \infty$.  In  the case $p=\infty$, the proof is similar.~~$\Diamond$


\begin{lem}\label{lemKK2-}
Let $1\le p\le \infty$, $1/p+1/q=1$, and $\{\w\phi_j\}_j \in \mathcal{L}_{q}$. Then, for any $f\in L_p$ and $j\in \N$, we have
\begin{equation*}
 \left\Vert \left\{\langle f, \w\phi_j(\cdot - M^{-j} k)\rangle \right\}_{k} \right\Vert_{\ell_{p,M^j}}\le \Vert \w\phi_j \Vert_{\mathcal{L}_{q,j}}\Vert f\Vert_p.
\end{equation*}
\end{lem}
{\bf Proof.}
In the case $p=\infty$, the proof is obvious since $\Vert \w\phi_j \Vert_{\mathcal{L}_{1,j}}=\Vert \w\phi_j \Vert_1$. For $p<\infty$, we have
\begin{equation*}
\begin{split}
  \bigg(\sum_{k\in D(M^j)} &|\langle f, \w\phi_j(\cdot - M^{-j} k)\rangle|^p\bigg)^\frac1p\\
  &=\(\sum_{k\in D(M^j)} \bigg|\sum_{l\in D(M^j)}
  \int_{\frac{l_1}{m_1^{j}}}^{\frac{l_1+1}{m_1^{j}}}dx_1\dots \int_{\frac{l_d}{m_d^{j}}}^{\frac{l_d+1}{m_d^{j}}}  f(x)\overline{\w\phi_j}(x-M^{-j}k) dx_d  \bigg|^p\)^\frac1p\\
  &=\(\sum_{k\in D(M^j)} \bigg|\int_{M^{-j}\T^d} \sum_{l\in D(M^j)}
    f(x+M^{-j}l)\overline{\w\phi_j}(x-M^{-j}(k-l)) dx  \bigg|^p\)^\frac1p\\
  &\le \int_{M^{-j}\T^d} \(\sum_{k\in D(M^j)} \bigg| \sum_{l\in D(M^j)}
    f(x+M^{-j}l)\overline{\w\phi_j}(x-M^{-j}(k-l))  \bigg|^p\)^\frac1p dx,
\end{split}
\end{equation*}
where the last formula follows from Minkowski's inequality.

Next, applying Young's inequality for the discrete convolution and H\"older's inequality, we derive that the last expression can be estimated from above by

\begin{equation*}
\begin{split}
  &\int_{M^{-j}\T^d}  \(\sum_{l\in D(M^j)} |f(x+M^{-j}l)|^p\)^\frac1p \sum_{k\in D(M^j)} |\w\phi_j(x-M^{-j}k)| dx\\
  &\le \Vert f\Vert_p \(\int_{M^{-j}\T^d}  \(\sum_{k\in D(M^j)} |\w\phi_j(x-M^{-j}k)|\)^q dx\)^\frac1q=m^{j/p}\Vert \w\phi_j \Vert_{\mathcal{L}_{q,j}}\Vert f\Vert_p,
\end{split}
\end{equation*}
which proves the lemma.~$\Diamond$

Now, we are ready to formulate and prove the main results of this section. We start from the case of strictly compatible functions/distributions $\phi_j$ and $\w\phi_j$.

\begin{theo}\label{cor1}
Let $1\le p\le\infty$, $1/p+1/q=1$, and $N\ge 0$. Suppose that $\{\w\phi_j\}_j\in \mathcal{S}_{N,p}'$,
$\{\phi_j\}_j\in \mathcal{B}$ with respect to the parameter $\delta \in (0,1/2)$,
and
$\phi_j$ and $\w\phi_j$  are strictly compatible with respect to $\delta$.
Then, for any $f\in  \mathbb{B}_{p,1}^{d/p+N}(M)$ and $j\in \N$, we have
\begin{equation}\label{KS000}
  \Vert f - Q_j(f,\phi_j,\w\phi_j) \Vert_p\le C K_{\phi_j,q} m^{-j(\frac1p+\frac Nd)}\sum_{\nu=j}^\infty m^{(\frac1p+\frac Nd)\nu} E_{\delta M^\nu}(f)_p,
\end{equation}
if, additionally, $\{\w\phi_j\}_j\in \mathcal{L}_{q}$, then for any $f\in L_p$, we have
\begin{equation}\label{KS000+++}
  \Vert f - Q_j(f,\phi_j,\w\phi_j) \Vert_p\le C K_{\phi_j,q} \Vert \w\phi_j \Vert_{\mathcal{L}_{q,j}} E_{\delta M^j}(f)_p,
\end{equation}
where the constant $C$ does not depend on $f$ and $j$.
\end{theo}

{\bf Proof.}
Let $T_j\in \mathcal{T}_{\delta M^j}$ be such that $\Vert f-T_j\Vert_p=E_{\delta M^j}(f)_p$, then
\begin{equation}\label{I123}
  \begin{split}
      &\bigg\Vert f - \frac1{m^j}\sum\limits_{k\in D(M^{j})} \langle f, \w\phi_j(\cdot - M^{-j} k) \rangle \phi_j (\cdot - M^{-j}k)\bigg\Vert_p\\
      &\le \Vert f-T_j\Vert_p+\bigg\Vert T_j-\frac1{m^j}\sum\limits_{k\in D(M^{j})} \langle T_j, \w\phi_j(\cdot - M^{-j} k) \rangle \phi_j (\cdot - M^{-j}k)\bigg\Vert_p\\
      &\qquad\qquad\quad\,\,\,\,+\bigg\Vert \frac1{m^j}\sum\limits_{k\in D(M^{j})} \langle f-T_j, \w\phi_j(\cdot - M^{-j} k) \rangle \phi_j (\cdot - M^{-j}k)\bigg\Vert_p:=I_1+I_2+I_3.
   \end{split}
\end{equation}


To estimate $I_2$, we note that by Theorem~\ref{theo_my1},
\begin{equation*}
  \begin{split}
      T_j(x)-\frac1{m^j}\sum\limits_{k\in D(M^{j})} \langle T_j, \w\phi_j(\cdot - M^{-j} k) \rangle \phi_j (x - M^{-j}k)=0,
   \end{split}
\end{equation*}
which implies that
\begin{equation}\label{I2}
  I_2=0.
\end{equation}

Consider $I_3$. Using Lemmas~\ref{lemKK1} and~\ref{lem1Be}, we derive
\begin{equation}\label{I3}
  \begin{split}
      I_3&\le C_1 K_{\phi_j,q}\Vert \{\langle f-T_j, \w\phi_{j}(\cdot-M^{-j}k)\rangle \}_{k} \Vert_{\ell_{p,M^j}}\\
	&\le C_2 K_{\phi_j,q} \sum_{\mu=n}^\infty \Vert \{\langle T_{\mu+1}-T_\mu, \w\phi_{j}(\cdot-M^{-j}k)\rangle \}_{k} \Vert_{\ell_{p,M^j}}\le
	 C_3 K_{\phi_j,q} \sum_{\mu=n}^\infty m^{\mu(\frac{N}d+\frac1p)}E_{\delta M^\mu}(f)_p.
   \end{split}
\end{equation}

Then, combining~\eqref{I123}, \eqref{I2}, and~\eqref{I3}, we prove~\eqref{KS000}.

To obtain inequality~\eqref{KS000+++}, it is sufficient to use inequalities~\eqref{I123} and~\eqref{I2} as well as the following estimate
\begin{equation}\label{dopdopdop}
  \begin{split}
      I_3&\le C_4 K_{\phi_j,q}\Vert \{\langle f-T_j, \w\phi_{j}(\cdot-M^{-j}k)\rangle \}_{k} \Vert_{\ell_{p,M^j}}\le C_5 K_{\phi_j,q} \Vert \w\phi_j \Vert_{\mathcal{L}_{q,j}}\Vert f-T_j\Vert_p,
   \end{split}
\end{equation}
which easily follows from  Lemmas~\ref{lemKK1} and~\ref{lemKK2-}.~~$\Diamond$

\medskip
Applying Hausdorff--Young's inequality~\eqref{fT12} to the right-hand sides of~\eqref{KS000} and~\eqref{KS000+++}, we derive the following improvements of the error estimate given in Theorem~\ref{theo_my1}.

\begin{coro}
Let $2\le p\le\infty$, $1/p+1/q=1$, and $N\ge 0$. Suppose that $M=\lambda I_d$, $\lambda>1$, $\{\w\phi_j\}_j\in \mathcal{S}_{N,p}'$,
$\{\phi_j\}_j\in \mathcal{B}$ with respect to the parameter $\delta \in (0,1/2)$,
and
$\phi_j$ and $\w\phi_j$  are strictly compatible with respect to $\delta$.
Then, for any $f\in  \mathbb{B}_{p,1}^{d/p+N}(M)$ and $j\in \N$, we have
\begin{equation*}
  \begin{split}
\Vert f - Q_j(f,\phi_j,\w\phi_j) \Vert_p\le C K_{\phi_j,q} \lambda^{-j(\frac dp+N)} \sum_{\nu=j}^\infty \lambda^{(\frac dp+N)\nu}\(\sum_{|r|\ge \delta\lambda^\nu} |\h f(r)|^q\)^{1/q}\\
   \end{split}
\end{equation*}
if, additionally, $\{\w\phi_j\}_j\in \mathcal{L}_{q}$, then
\begin{equation*}
  \begin{split}
\Vert f - Q_j(f,\phi_j,\w\phi_j) \Vert_p\le C K_{\phi_j,q}\Vert \w\phi_j \Vert_{\mathcal{L}_{q,j}} \(\sum_{|r|\ge \delta\lambda^j} |\h f(r)|^q\)^{1/q},\\
   \end{split}
\end{equation*}
where the constant $C$ does not depend on $f$ and $j$.
\end{coro}

\bigskip

In light of Proposition~\ref{coro1}, it is not difficult to derive the following improvements of Theorem~\ref{cor1}, in which we replace the best approximation $E_{\delta M^j}(f)_p$ by
$$
E_{M^j}^*(f)_p:=\inf\left\{\Vert f-T\Vert_p\,:\, \spec T \subset D(M^j)\right\},\quad j\in \N.
$$

\begin{prop}\label{propLag}
If under the assumptions of Theorem~\ref{cor1}, equality~\eqref{my3} holds for all $l \in D(M^j)$ and
$\spec \vp_j \subset D(M^j)$, then, for any $f\in \mathbb{B}_{p,1}^{d/p+N}(M)$ and $j\in \N$, we have
\begin{equation*}
  \Vert f - Q_j(f,\phi_j,\w\phi_j) \Vert_p\le C K_{\phi_j,q} m^{-j(\frac1p+\frac Nd)}\sum_{\nu=j}^\infty m^{(\frac1p+\frac Nd)\nu} E_{M^\nu}^*(f)_p
\end{equation*}
if, additionally, $\{\w\phi_j\}_j\in \mathcal{L}_{q}$, then for any $f\in L_p$, we have
\begin{equation*}
  \Vert f - Q_j(f,\phi_j,\w\phi_j) \Vert_p\le C K_{\phi_j,q} \Vert \w\phi_j \Vert_{\mathcal{L}_{q,j}} E_{M^j}^*(f)_p,
\end{equation*}
where the constant $C$ does not depend on $f$ and $j$.
\end{prop}

Now, we consider the case of weakly compatible functions/distributions $\phi_j$ and $\w\phi_j$.

\begin{theo}\label{cor3}

Let $2\le p\le \infty$, $1/p+1/q=1$, and $N\ge 0$. Suppose that $\{\w\phi_j\}_j\in \mathcal{S}_{N,p}'$, $\{\phi_j\}_j\in \mathcal{B}$ with respect to the parameter $\delta \in (0,1/2)$,
and  $\phi_j$ and $\w\phi_j$  are weakly compatible of order $s>0$.
Then, for any $f\in \mathbb{B}_{p,1}^{d/p+N}(M)$ and $j\in \N$, we have
\begin{equation}\label{-KS000+}
\begin{split}
     \Vert f - &Q_j(f,\phi_j,\w\phi_j) \Vert_p\le C\(\Vert M^{-j}\Vert^s  \|f\|^{In}_{A_q^{s},j} +K_{\phi_j,q} m^{-j(\frac1p+\frac Nd)}\sum_{\nu=j}^\infty m^{(\frac1p+\frac Nd)\nu} E_{\delta M^\nu}(f)_p\)
\end{split}
\end{equation}
if, additionally, $\{\w\phi_j\}_j\in \mathcal{L}_{q}$, then for any $f\in L_p$, we have
\begin{equation}\label{--KS000+}
\begin{split}
     \Vert f - Q_j(f,\phi_j,\w\phi_j) \Vert_p\le C\(\Vert M^{-j}\Vert^s\|f\|^{In}_{A^{s},j} + K_{\phi_j,q}\Vert \w\phi_j \Vert_{\mathcal{L}_{q,j}} E_{\delta M^j}(f)_p\),
\end{split}
\end{equation}
where the constant $C$ does not depend on $f$ and $j$.
\end{theo}

{\bf Proof.} First, we prove inequality~\eqref{-KS000+}. Consider the de la Vall\'ee--Poussin means $V_{j}(f)$ defined by
$$
V_j(f)(x)=\sum_{k\in \Z^d} v_\delta(M^{-j}k)\h f(k)e^{2\pi i(k,x)},
$$
where $v_\delta \in C^\infty(\R^d)$, $v_\delta(\xi)=1$ if $|\xi|<\delta$ and $v_\delta(\xi)=0$ if $\xi \not\in (-\frac12,\frac12)^d$.
As usual, we have for any $j\in \N$ that
\begin{equation}\label{finVal}
  \Vert f-V_j(f)\Vert_p\le (1+\Vert V_j\Vert_1)E_{\delta M^j}(f)_p\le c E_{\delta M^j}(f)_p,
\end{equation}
where the constant $c$ does not depend on $j$ and $f$.

Repeating the proof of Theorem~\ref{cor1} with $V_j(f)$ instead of the polynomials of the best approximation $T_j$ and using~\eqref{finVal}, we derive
\begin{equation}\label{Vpart1}
  \begin{split}
      &\bigg\Vert f - \frac1{m^j}\sum\limits_{k\in D(M^{j})} \langle f, \w\phi_j(\cdot - M^{-j} k) \rangle \phi_j (\cdot - M^{-j}k)\bigg\Vert_p\\
      &\le \bigg\Vert V_j(f)-\frac1{m^j}\sum\limits_{k\in D(M^{j})} \langle V_j(f), \w\phi_j(\cdot - M^{-j} k) \rangle \phi_j (\cdot - M^{-j}k)\bigg\Vert_p\\
      &\qquad\qquad\qquad\qquad+C_1 K_{\phi_j,q} m^{-j(\frac1p+\frac Nd)}\sum_{\nu=j}^\infty m^{(\frac1p+\frac Nd)\nu} E_{\delta M^\nu}(f)_p:=J_1+J_2.
  \end{split}
\end{equation}
By Theorem~\ref{big_p_th}, taking into account that the Strang-Fix conditions of order $s$ for $\phi_j$
and the weak compatibility conditions for $\phi_j$ and $\w \phi_j$ of order $s$ are satisfied (see the proof of Theorem~\ref{theo_my1}),  we derive
\begin{equation}\label{Vpart2}
  \begin{split}
J_1 \le C_2\|M^{-j}\|^{s}\|V_j(f)\|^{In}_{A_q^{s},j}\le C_3\|M^{-j}\|^{s}\|f\|^{In}_{A_q^{s},j}.
  \end{split}
\end{equation}
Thus, combining~\eqref{Vpart1} and~\eqref{Vpart2}, we prove~\eqref{-KS000+}.

The proof of estimate~\eqref{--KS000+} is similar. One needs only to use inequality~\eqref{dopdopdop} instead of~\eqref{I3}.~~$\Diamond$

\medskip

Next, using the Hausdorff-Young inequality, we obtain the following corollaries, which provide two shaper versions of Theorem~\ref{big_p_th}.

\begin{coro}
If under the assumptions of Theorem~\ref{cor3}, $\{\w\phi_j\}_j\in \mathcal{L}_{q}$ and $f\in A_q^\gamma$, where $\gamma>0$, then 
$$
\Vert f - Q_j(f,\phi_j,\w\phi_j) \Vert_p\le CK_{\phi_j,q} \Vert M^{-j}\Vert^{\min\{s,\gamma\}}\|f\|_{A_q^\gamma}.
$$
In particular, if $M=\lambda I_d$, $\sup_j K_{\phi_j,q} <\infty$, and $\h f(n)=\mathcal{O}(|n|^{-\kappa})$ for some $\kappa>d/q$, then
$$
\Vert f - Q_j(f,\phi_j,\w\phi_j) \Vert_p= \mathcal{O}( \lambda^{-j\min\{s,\kappa-d/q\}}).
$$
\end{coro}

\begin{coro}\label{corsharp}
If under the assumptions of Theorem~\ref{cor3}, $M=\lambda I_d$, $0<s\le d/p+N$, $f\in A_q^s \cap \mathbb{B}_{p,1}^{d/p+N}(M)$, then
$$
\Vert f - Q_j(f,\phi_j,\w\phi_j) \Vert_p\le CK_{\phi_j,q}  \lambda^{-js}.
$$
\end{coro}


\begin{rem}
Note that this result provides a shaper version of  Theorem~\ref{big_p_th} because  there exist functions $f\not\in A_q^\gamma$, $\gamma>d/p+N$, for which conditions of Corollary~\ref{corsharp} are valid. As an example in the case $d=1$, one can take the function
$$
f(x)=\sum_{k=1}^\infty \frac{e^{ik\log k}}{k^{\frac12+\frac1p+N+\varepsilon}} e^{2\pi i kx},\quad \varepsilon>0.
$$
It follows from~\cite[Ch. V, Theorem 4.2]{Z} that $f^{(N)}\in {\rm Lip}(\frac1p+\varepsilon)$ and, therefore, by the classical Jackson inequality (see, e.g.,~\cite[p.~260]{Timan}),  we have that $E_{\lambda^\nu}(f)_p=\mathcal{O}(\lambda^{-\gamma \nu})$ with $\gamma=1/p+N+\varepsilon$, which implies that $f\in B_{p,1}^{1/p+N}$. At the same time, $f\in A_q^s$ for $s<\gamma-1/q$
and $f\not\in A_q^\gamma$.
\end{rem}

\section{Examples}

1. In fact, for a sequence of periodic distributions $\{\w\phi_j\}$ satisfying conditions~(\ref{fWPhi}) and such that $|\h{\w \phi_j}(k)| \ge c > 0$ for $k\in\zd: |M^{*-j} k| \le \delta$ for any $j\in\n$ and some fixed $\delta \in (0,\frac 12)$, an appropriate strictly compatible sequence of trigonometric polynomials $\phi_j$ can be constructed via defining its Fourier coefficients by condition~(\ref{my3}), namely
$$
\overline{\h{\phi_j} (k)} = \frac{1}{\h{\w \phi_j}(k)}, \quad k\in\zd: |M^{*-j} k| \le \delta
$$
and $\h{\phi_j} (k) = 0$ for others $k \in \zd$. Obviously, the sequence  $\{\phi_j\}_j$ belongs to $ \mathcal{B}$, and hence
we are under the assumptions of Theorem~\ref{theo_my1} or Theorem~\ref{cor1}
(if $\{\w\phi_j\}_j \in \mathcal{S}_{N,p}'$ and  $M$ is diagonal).
For instance, assume that $\w\phi_j$ is a periodic distribution
corresponding to some differential operator. Namely, let
$$
\h{\w\phi_j} (k) = \sum_{[\beta]\le N} c_\beta (2\pi i M^{*-j} k)^{\beta}, \quad N\in \z_+,
\quad c_\nul  \neq 0.
$$
For any good enough function $f$, we have
\begin{equation*}
  \begin{split}
       \langle f,\w\phi_j(\cdot - M^{-j}k)\rangle &= \sum_{l\in\zd} \h f(l) \overline{\h{\w\phi_j} (l)} e^{2\pi i (M^{-j}k,l)} =
      \sum_{[\beta]\le N} \overline{c}_\beta \sum_{l\in\zd} \h f(l) (-2\pi i M^{*-j} l)^{\beta} e^{2\pi i (k,M^{*-j}l)} \\
&=\sum_{[\beta]\le N} \overline{c}_\beta  [D^{\beta} f(-M^{-j} \cdot)] (k) =:  [L f(M^{-j} \cdot) ] (k).
   \end{split}
\end{equation*}
This sequence  of periodic distributions  $\w\phi_j$ satisfies conditions~(\ref{fWPhi}).
Since $\h{\w\phi_j}(k)=P(M^{*-j} k)$, where $P$
is an algebraic polynomial and $P(\nul)\ne0$, there exists $\delta>0$ such that
$|\h{\w\phi_j} (k)|  \ge c > 0$ whenever $k\in\zd: |M^{*-j} k| \le \delta.$
In  particular, if $N=0$ and $c_\nul = 1$, then $\h{\w\phi_j} (k) = 1$ for $k\in\zd$.
For the strict compatibility, one can take $\h{\phi_j} (k) = 1$
for $k\in\zd: |M^{*-j} k| \le \frac12$ and $\h{\phi_j} (k) = 0$ for other $k \in \zd$,
i.e., $\{\phi_j\}_j$ is a sequence of Dirichlet-type kernels.

2. Let $\h{\w\phi_j} (k) = 1$ for $k\in\zd$. In order to achieve only the  weak compatibility, we can
take truncated Fejer-type kernels $\phi_j$ defined by
$$
\h{\phi_j}(k) =\left\{
                 \begin{array}{ll}
                   1 - C_F \|M^{*-j}k\|_{\infty}, & \hbox{if $|M^{*-j} k| \le \delta$,} \\
                   0, & \hbox{otherwise,}
                 \end{array}
               \right.
$$
where $\delta \in (0,\frac 12)$ and $C_F$ is a positive real number. In this case
$$
1 -  \h{\w\phi_j} (k) \h{\phi_j}(k) =  C_F \|M^{*-j}k\|_{\infty} \le  C_F |M^{-j}k|,  \quad k\in D(M^{*j}),
$$
which means that compatibility condition~(\ref{fPhiPhiW}) is valid for $s=1.$ Thus, we are in case of Theorem~\ref{cor3}.

Alongside, consider the following Fejer-type kernels defined by
$$
\h{\phi_j}(k) =\left\{
                 \begin{array}{ll}
                   1 - C_F \|M^{*-j}k\|_{\infty}, & \hbox{if $\|M^{*-j}k\|_{\infty} \le 1$,} \\
                   0, & \hbox{otherwise.}
                 \end{array}
               \right.
$$
In this case, we have to check the Strang-Fix conditions only for points $M^{*j}n + r$, where $\|n\|_{\infty} = 1$ and $r\in D(M^{*j})$. Then
$$
| \h{\phi_j}(M^{*j}n + r)| = |1 - \|n + M^{*-j}r\|_{\infty}| = |\|n\|_{\infty} - \|n + M^{*-j}r\|_{\infty}| \le \|M^{*-j}r\|_{\infty} \le |M^{-j}r|.
$$
Therefore, the Strang-Fix condition of order $1$ for the sequence $\{\phi_j\}_j$ is satisfied. So, we are under assumptions of Theorem~\ref{theoMain} or Theorem~\ref{big_p_th} with $s=1$.

3.  Next, we discuss  sequences $\{\phi_j\}_j$ obtained  by periodization of  splines, which are applicable in Theorem~\ref{theoMain} or~\ref{big_p_th}. We restrict ourselves to the case $d=1$, and $M$ is an integer greater than 1.
Assume that  $\phi^{s}$ is a B-spline of an even order $s$, whose Fourier transform is given by
\begin{equation}
\h{\phi^{s}} (\xi) =  \left( \frac{\sin \pi \xi}{\pi \xi}\right)^{s}.
\label{fBSpline}
\end{equation}
Recall that ${\rm supp\,}\phi^{s}\subset [-s/2,s/2]$
and the non-periodic Strang-Fix conditions of order $s$ are valid for $\phi^{s}$.
Therefore, after periodization of $\phi^s$ by~(\ref{fPeriodiz}), we obtain the functions $\phi^{s}_{j}(x)$
whose Fourier coefficients are
$$
\h{\phi^{s}_{j}} (k) = \left( \frac{\sin \pi \frac{k}{M^j}}{\pi \frac{k}{M^j}}\right)^{s}.
$$
In fact, after periodization for a big enough $j$ (such that $[-\frac{s/2}{M^j},\frac{s/2}{M^j}] \in [-1/2,1/2]$) we have
 $\phi^{s}_{j}(x) = \phi^{s}(M^{j}x)$ for $x\in[-1/2,1/2]$, i.e. $\phi^{s}_{j}$ is a contracted B-spline  $\phi^s$.

Now, check that the periodic Strang-Fix conditions of  order $s$  for the sequence $\{\phi^{s}_{j}\}_j$  are valid and the corresponding sequence $\{b_n\}_n$ from~(\ref{fPhiSF}) belongs to $\ell_q$ for any $q>1$.
Indeed, let $n\neq 0$ and $r\in D(M^j)$, then
$$
| \h{\phi^{s}_j}(M^{j}n+r)| =  \left| \frac{\sin \pi \frac{M^{j}n+r}{M^j}}{\pi \frac{M^{j}n+r}{M^j}}\right|^{s} =  \left| \frac{\sin \pi \frac{r}{M^j}}{\pi (n +  \frac{r}{M^j})}\right|^{s} \le  \left| \frac{\frac{r}{M^j}}{n +  \frac{r}{M^j}}\right|^{s} \le
\frac{2^s}{(2|n|-1)^s}\left| \frac{r}{M^j} \right|^s.
$$
Also, we have that $b_n = \frac{2^s}{(2|n|-1)^s}$, $n\neq 0$, which implies that $\{b_n\}_n \in \ell_q$.

Next, for a given sequence $\{\w\phi_j\}_j$, we can consider
linear combinations of splines in order to construct an appropriate sequence  $\{\phi_j\}_j$. For instance,  in order to reduce noise contribution,  it is reasonable to use smoothed version of samples instead of the exact samples of $f$ (see, e.g.~\cite{ZWS}). Let $\w\phi_j$  be such that
$$
\langle f,\w\phi_j(\cdot - M^{-j}k)\rangle = \frac 12 f(M^{-j}k) + \frac 14 f(M^{-j}(k+1)) + \frac 14 f(M^{-j}(k-1)).
$$
This means that $\h{\w\phi_j}(k) =\frac 12 + \frac 12 \cos 2\pi \frac{k}{M^j}.$  Thus,  condition~(\ref{fWPhi}) is satisfied with $N=0.$ Consider a periodized B-spline as a dual sequence $\{\phi^s_j\}$, whose Fourier coefficients are defined by~(\ref{fBSpline}).
Therefore,
$$
1- \h{\w\phi_j} (k) \h{\phi^s_j}(k) = 1- \left(\frac 12 + \frac 12 \cos 2\pi \frac{k}{M^j}\right) \left( \frac{\sin \pi \frac{k}{M^j}}{\pi \frac{k}{M^j}}\right)^{s}.
$$
Since, by Taylor's formula, $f_1(x) := 1-\left(\frac 12 + \frac 12 \cos 2\pi x\right) \left( \frac{\sin \pi x}{\pi x}\right)^{s} = \frac 16 \left(6 \pi^2+s \pi ^2\right) x^2 + \mathcal{O}(x^4)$ as $x\to 0$, it is clear that the order of compatibility cannot be better than 2. Thus, using Lagrange's reminder for the case $s=2$, we can state that
$$
|f_1(x)| \le \frac{|x|^2}{2!} \max_{x\in[-1/2,1/2]} |f_1''(x)|  = \frac 43 \pi^2 |x|^2,
$$
which implies
$$
|1- \h{\w\phi_j} (k) \h{\phi^s_j}(k)| \le \frac 43 \pi^2 \left|\frac{k}{M^j}\right|^2, \quad k \in D(M^{*j}).
$$

Now, consider a linear combination of the shifted splines $\phi_j^{s}$, which allows to provide a better order of compatibility. For instance, let $\phi_j := u_1 \phi_j^{s} + \frac{u_2}{2} \phi_j^{s}(\cdot + \frac{1}{M^{j}}) + \frac{u_2}{2} \phi_j^{s}(\cdot - \frac{1}{M^{j}}).$ The Fourier coefficients of $\phi_j$ are
\be
\label{30}
\h{\phi_{j}} (k) =  \left(  u_1\left(\frac{\sin \pi \frac{k}{M^j}}{\pi \frac{k}{M^j}}\right)^{s} + u_2 \cos\Big(2 \pi \frac{k}{M^j}\Big) \left(\frac{\sin \pi \frac{k}{M^j}}{\pi \frac{k}{M^j}}\right)^{s} \right).
\ee
To find the compatibility order, we consider the function
$$
f_2(x) :=  1-\left(\frac 12 + \frac 12 \cos 2\pi x\right) \left( \frac{\sin \pi x}{\pi x}\right)^{s} (u_1 + u_2 \cos (2 \pi x)).
$$
From its  Taylor's formula near the origin, we get
$$
f_2(x) = (1-u_1-u_2) + \frac{\pi^2}{6} x^2 \left( u_1 s+6 u_1+ u_2 s+18  u_2 \right) + \mathcal{O}(x^4).
$$
It is clear that for $s=4$, $u_1 = 11/6$, $u_2 = -5/6$, the first two terms are vanished.
In order to get constant value $b_0$ from condition~(\ref{fPhiPhiW}),
Moreover, using Taylor's reminder, we have
$$
|f_2(x)| \le \frac{|x|^4}{4!} \max_{x\in[-1/2,1/2]} |f_2^{(IV)}(x)|  = \frac{32}{15} \pi^4 |x|^4.
$$
This yields the weak compatibility of order 4 for  $\{\phi_j^s\}_j$ and $\{\w\phi_j\}_j$, i.e.,
$$
|1- \h{\w\phi_j} (k) \h{\phi_j}(k)| \le \frac{32}{15} \pi^4 \left|\frac{k}{M^j}\right|^4, \quad k \in D(M^{*j}).
$$

A similar procedure can be applied for the construction of a sequence  $\{\phi_j\}_j$, which is weakly compatible with a sequence of distributions $\{\w\phi_j\}_j$ corresponding to some differential operator. For instance, assume that
$$
\h{\w\phi_j} (k) = 1+ c_1 (2\pi i  M^{-j} k)^{2}
$$
and
$\h{\phi_j}(k)$
are defined by~(\ref{30}) with $s=4$.
Let us try to achieve  the weak compatibility of order 4 for the corresponding sequences  $\{\phi_j\}_j$ and $\{\w\phi_j\}_j$.
For this, it is enough to apply Taylor's formula to the function 
$$
f_3(x) :=  1-\left(1+ c_1 (2\pi i x)^{2}\right) \left( \frac{\sin \pi x}{\pi x}\right)^{4} (u_1 + u_2 \cos (2 \pi x))
$$
and chose $u_1 = 2c_1 + 4/3$ and $u_2 = 1-u_1$ in (\ref{30}). In particular, if $c_1 = -1/4,$ then
$$
|1-\h{\w\phi_j} (k) \h{\phi_j}(k)| \le \frac{7}{15} \pi^4 \left|\frac{k}{M^j}\right|^4, \quad k \in D(M^{*j}),
$$
which implies the required fact. 

\end{document}